\documentclass[11pt,a4paper]{article}
\usepackage{amsmath,bm,graphicx,booktabs}
\usepackage{geometry,indentfirst,float,appendix,lmodern}
\usepackage{multirow}
\usepackage[authoryear]{natbib}
\usepackage{algorithm,algorithmic}
\usepackage{hyperref}
\usepackage[affil-it]{authblk}
\usepackage{setspace}
\usepackage{inputenc}
\usepackage[english]{babel}
\usepackage{amsthm}
\usepackage{color}
\usepackage{courier}
\usepackage{lmodern}
\usepackage{amssymb}
\usepackage{mathrsfs}
\usepackage{enumitem}

\providecommand{\keywords}[1]{\textbf{\textit{key words: }} #1}

\newtheorem{theorem}{Theorem}

\newtheorem{remark}{Remark}
\newtheorem{corollary}{Corollary}
\newtheorem{lemma}{Lemma}

\newcommand{\pcite}[1]{\citeauthor{#1}'s \citeyearpar{#1}}

  \title{Convergence analysis of the block Gibbs sampler for Bayesian
    probit linear mixed models with improper priors} 
    
  \date{}

\author[1]{Xin Wang \thanks{Email: wangx172@miamioh.edu}}
\author[2] {Vivekananda Roy \thanks{Email:vroy@iastate.edu }}
\affil[1] {Department of Statistics, Miami University}
\affil[2]{Department of Statistics, Iowa State University}

\begin{document}

\maketitle

\begin{abstract}
  In this article, we consider Markov chain Monte Carlo (MCMC)
  algorithms for exploring the intractable posterior density
  associated with Bayesian probit linear mixed models under improper
  priors on the regression coefficients and variance components. In
  particular, we construct a two-block Gibbs sampler using the data
  augmentation (DA) techniques. Furthermore, we prove geometric
  ergodicity of the Gibbs sampler, which is the foundation for
  building central limit theorems for MCMC based estimators and
  subsequent inferences. The conditions for geometric convergence are
  similar to those guaranteeing posterior propriety.  We also provide
  conditions for the propriety of posterior distributions with a general link function when the design matrices take
  commonly observed forms. In general, the Haar parameter expansion
  for DA (PX-DA) algorithm is an improvement of the DA algorithm and
  it has been shown that it is theoretically at least as good as the
  DA algorithm. Here we construct a Haar PX-DA algorithm, which has
  essentially the same computational cost as the two-block Gibbs
  sampler.
\end{abstract}

\keywords{Data augmentation, Drift condition, Geometric ergodicity, GLMM, Haar PX-DA algorithm, Markov chains, Posterior propriety}

\section{Introduction}
\label{sec_introduction}

Generalized linear mixed models (GLMMs) are generalized linear models
with random terms in the linear predictor. The random effects in the
GLMM can accommodate for overdispersion often present in non-Gaussian
data, and dependence among correlated observations arising from
longitudinal or repeated measures studies. GLMM is one of the most
frequently used statistical models. Here, we consider a popular Bayesian
GLMM for binary data, namely, the probit linear mixed model.

Let $(Y_1, Y_2, \dots, Y_n)$ denote the vector of Bernoulli
random variables. Let $\bm{x}_i$ and $\bm{z}_i$ be the
$p \times 1$ and $q \times 1$ known covariates and random effect
design vectors respectively associated with the $i$th observation for $i=1,\dots,n$.
Let $\bm{\beta} \in \mathbb{R}^p$ be the unknown vector of regression
coefficients and $\bm{u} \in \mathbb{R}^q$ be the random effects
vector.  A GLMM can be built \citep{mcculloch2011model, breslow1993approximate} with a link
function that connects the expectation of $Y_i$ with 
$\bm{x}_i$ and $\bm{z}_i$. One of the very popular link functions is the probit
link function, $\Phi^{-1}$, resulting in  
\begin{equation}
\label{eq_link}
P(Y_i =1) = \Phi(\bm{x}_i^T\bm{\beta} + \bm{z}_i^T\bm{u}),
\end{equation}
 where $\Phi(\cdot)$ is the cumulative distribution function of the
standard normal random variable. Assume that we have $r$ random
effects with $\bm{u} = (\bm{u}_1^T,\dots, \bm{u}_r^T)^T$, where
$\bm{u}_j$ is a $q_j\times1$ vector with $q_j >0$, $q_1+\cdots + q_r = q$, and
$\bm{u}_j \stackrel{\text{ind}}{\sim} N(0, \bm{I}_{q_j}1/\tau_j)$,
where $\tau_j\in \mathbb{R}_+ \equiv (0,\infty)$ is the precision parameter
associated with $\bm{u}_j$ for $j=1,\dots,r$.  Let
$\bm{\tau}= (\tau_1, \dots, \tau_r)$, thus, the data model for the
probit GLMM is
\begin{eqnarray}
\label{eq:data}
Y_{i} | \bm{\beta}, \bm{u} , \bm{\tau}& \overset{\text{ind}}\sim & \text{Bern}(\alpha_{i}) \text{ for } i=1,\dots,n \;\text{ with} \nonumber\\
\alpha_{i} & = & \Phi(\bm{x}_{i}^{T}\bm{\beta}+\bm{z}_{i}^{T}\bm{u}) \text{ for }i=1,\dots,n, \\
\bm{u}_{j}| \bm{\beta}, \bm{\tau} & \overset{\text{ind}}\sim & N\Big(0,\frac{1}{\tau_{j}}\bm{I}_{q_{j}}\Big) , \, j=1,\dots,r .\nonumber
\end{eqnarray}
 Let
$\bm{y} = (y_1, y_2, \dots, y_n)^T$ be the observed Bernoulli response variables. Note that,
the likelihood function for $(\bm{\beta}, \bm{\tau})$ is
\begin{eqnarray}
  \label{eq:lik}
  L(\bm{\beta}, \bm{\tau} | \bm{y}) &=& \int_{\mathbb{R}^{q}} \prod_{i=1}^n\left[\Phi(\bm{x}_i^T\bm{\beta} + \bm{z}_i^T\bm{u})\right]^{y_i}\left[1-\Phi(\bm{x}_i^T\bm{\beta} +\bm{z}_i^T\bm{u})\right]^{1-y_i} \\
  & \times & \phi_q(\bm{u}; \bm{0}, \bm{D}(\bm{\tau})^{-1})d\bm{u}, \nonumber
\end{eqnarray}
which is not available in closed form. 
Here, $\phi_q(s; a, B)$ denotes the probability density
function of the $q-$dimensional normal distribution with mean vector
$a$, covariance matrix $B$ and evaluated at $s$, and
$\bm{D}(\bm{\tau}) = \oplus_{j=1}^r \tau_j \bm{I}_{q_j}$.

In Bayesian framework, one needs to specify the prior distributions of
$\bm{\beta}$ and $\bm{\tau}$. Assume $\bm{\beta}$ and $\bm{\tau}$ are apriori independent.  Let $\pi(\bm{\beta}) $ and $\pi(\bm{\tau})$ be the prior densities of $\bm{\beta}$ and $\bm{\tau}$ respectively.  Thus, the joint posterior
density of $(\bm{\beta},\bm{\tau})$ is
\begin{equation}
  \label{eq:postbetphi}
  \pi(\bm{\beta},\bm{\tau}|\bm{y}) = \frac{1}{c(\bm{y})} L(\bm{\beta}, \bm{\tau}| \bm{y})\pi(\bm{\beta})\pi({\bm{\tau}}),
\end{equation}
where
\begin{equation*}
c(\bm{y}) =  \int_{\mathbb{R}_+^{r}} \int_{\mathbb{R}^{p}} L(\bm{\beta}, \bm{\tau}| \bm{y}) \pi(\bm{\beta}) \pi(\bm{\tau})d\bm{\beta} d\bm{\tau},
\end{equation*}
is the marginal density of $\bm{y}$. Since the likelihood function
$L(\bm{\beta}, \bm{\tau}| \bm{y})$ is not available in closed form,
the posterior density is intractable for any choice of the prior
distributions of $\bm{\beta}$ and $\bm{\tau}$.  In this article, we consider an improper flat prior for $\bm{\beta}$, that is, $\pi(\bm{\beta})\propto 1$ and $\tau_j$'s, $j=1,\dots, r$, are apriori independent with
\begin{equation}
\label{eq:improper_prior}
\pi(\tau_j )\propto e^{-b_j\tau_j}\tau_j^{a_j-1},
\end{equation}
which can be proper or improper. In section~\ref{sec_propriety}, we
discuss conditions under which the posterior density
\eqref{eq:postbetphi} is proper, that is $c(\bm{y}) < \infty$. Generally, Markov chain Monte Carlo (MCMC)
algorithms are used for exploring the posterior density
\eqref{eq:postbetphi}.

Even in the absence of random effects, for the probit regression model,
the posterior distribution of $\bm{\beta}$ is difficult to sample from
\citep{roy2007convergence}. \pcite{albert1993bayesian} MCMC algorithm
for sampling from the posterior distribution associated with the
probit regression model is the most widely used data augmentation (DA)
algorithm. The DA technique used in \cite{albert1993bayesian}  can also be applied
to the probit linear mixed model. Following \cite{albert1993bayesian},
let $v_i \in \mathbb{R}$ be the continuous latent variable
corresponding to the $i$th binary observation $Y_i$, such that $Y_i = I(v_i >0)$, where
$v_i | \bm{\beta}, \bm{u}, \bm{\tau}\overset{\text{ind}} \sim N(\bm{x}_i^{T} \bm{\beta} +
\bm{z}_{i}^{T}\bm{u}$, 1) for $i =1,\dots, n$. Then
\begin{equation}
  \label{eq:vy}
  P(Y_i=1) = P(v_i >0) =\Phi(\bm{x}_i^{T} \bm{\beta} +
\bm{z}_{i}^{T}\bm{u}),
\end{equation}
that is, $Y_{i} | \bm{\beta}, \bm{u}, \bm{\tau} \overset{\text{ind}}\sim  \text{Bern}(\alpha_{i})$ as in \eqref{eq:data}. Note that
$\bm{v}| \bm{\beta}, \bm{u}, \bm{\tau} \sim N(\bm{X}\bm{\beta} + \bm{Z}\bm{u},\bm{I}_n)$, where $\bm{v} = (v_1,\dots,v_n)^T$, 
$\bm{X}_{n\times p} =(\bm{x}_1,\dots, \bm{x}_{n})^T $ and
$\bm{Z}_{n\times q} = (\bm{z}_1,\dots, \bm{z}_n)^{T}$.

Using the latent variables $\bm{v}$, we can introduce a joint density
$\pi(\bm{\beta}, \bm{u},\bm{v},\bm{\tau}|\bm{y})$ (see
section~\ref{sec_twoblock} for details) such that
\begin{equation}
  \label{eq:margjt}
  \int_{\mathbb{R}^{q}} \int_{\mathbb{R}^{n}} \pi(\bm{\beta}, \bm{u},\bm{v},\bm{\tau}|\bm{y}) d\bm{v}d\bm{u} = \pi(\bm{\beta},\bm{\tau}|\bm{y}),
\end{equation}
where $\pi(\bm{\beta},\bm{\tau}|\bm{y})$ is the posterior density
defined in \eqref{eq:postbetphi}. If all the full conditionals of the
joint density $\pi(\bm{\beta}, \bm{u},\bm{v},\bm{\tau}|\bm{y})$ are
easy to sample from, then a Gibbs sampler can be run and it can be
used to make inferences on the posterior density \eqref{eq:postbetphi}. Indeed this full
Gibbs sampler is traditionally used in the analysis of Bayesian probit
linear mixed models \citep{baragatti2011bayesian}. In this article,
instead of using full conditional distributions, we construct a 
two-block Gibbs sampler with $ \bm{\eta} \equiv (\bm{\beta}^T,\bm{u}^T)^T$ as one block and
$(\bm{v}^T, \bm{\tau}^T)^T$ as the other block --- which is our first contribution. In general, block Gibbs samplers are known to be better than the Gibbs samplers based on full conditional
distributions in terms of having smaller operator norm \citep{liu1994covariance}.

The above mentioned block Gibbs sampler has an everywhere strictly
positive Markov transition density, implying that the underlying
Markov chain is Harris ergodic \citep{asmu:glyn:2011,
  meyn1993markov}. Thus, the time average estimators based on the
block Gibbs sampler can be used to consistently estimate the
(posterior) means with respect to the joint density
$\pi(\bm{\beta}, \bm{u},\bm{v},\bm{\tau}|\bm{y})$. In practice, it is
crucial to know whether the Monte Carlo errors associated with these
estimates are sufficiently small. However, in order to provide valid
standard errors, we need to establish a central limit theorem (CLT)
for the time average estimators. Unlike for the ordinary Monte Carlo
methods based on iid samples, mere existence of the finite second
moment does not guarantee a CLT for MCMC estimators. One standard
method of establishing a CLT for MCMC estimators is to prove that the
underlying Markov chain is {\it geometrically ergodic}
\citep{jones2001honest}. Geometric ergodicity is also needed for
consistently estimating the asymptotic variance in the Markov chain
CLT \citep{flegal2010batch}.  \cite{roy2007convergence} and
\cite{chakraborty2017convergence} proved geometric ergodicity of \pcite{albert1993bayesian}
DA algorithm for the Bayesian probit
regression model under improper and proper priors on the regression
coefficients. For linear models, \cite{jones2004sufficient} and
\cite{tan2009block} analyzed the Gibbs sampler for one-way random
effects models under proper priors and improper priors
respectively. \cite{johnson2010gibbs} analyzed the block Gibbs sampler
for Bayesian linear mixed models under the assumption
$\bm{X}^T\bm{Z} = \bm{0} $. \cite{roman2012convergence} and
\cite{roman2015geometric} established geometric rate of convergence of
the Gibbs samplers for Bayesian linear mixed models under improper and
proper priors without the assumption of $\bm{X}^T\bm{Z} = \bm{0}
$.  Our second contribution, in this paper, is
establishing geometric convergence rates for the block Gibbs sampler
for Bayesian probit linear mixed models under improper priors.

DA algorithms are known to suffer from slow convergence
\citep{meng:vand:1999, vand:meng:2001}.  \cite{liu1999parameter}
proposed the {\it parameter expansion for data augmentation} (PX-DA)
algorithm, which can converge faster than the DA algorithm without
much extra computational effort \citep{vand:meng:2001,
  roy2014efficient}. \cite{hobert2008theoretical} proved that the Haar
PX-DA algorithm, that is based on a Haar measure, is better than any
other PX-DA algorithm and the original DA algorithm in both the
efficiency ordering and the operator norm ordering. For the probit
regression model, \cite{roy2007convergence}, through an example,
showed that the Haar PX-DA algorithm can lead to huge gains in
efficiency over the DA algorithm of \cite{albert1993bayesian}. Our
third contribution is to construct a Haar PX-DA algorithm improving
the block Gibbs sampler mentioned before.  Since geometric ergodicity
of the Haar PX-DA algorithm follows from geometric ergodicity of the
DA algorithm \citep{hobert2008theoretical}, we have CLTs for the Haar
PX-DA algorithm based estimators as well.

The article is organized as follows. In section \ref{sec_propriety},
we establish conditions for propriety of the posterior distribution \eqref{eq:postbetphi}
under improper priors, when $\bm{X}$ and $\bm{Z}$ take commonly
observed forms. The results in section \ref{sec_propriety} hold for a
general link function, not necessarily the probit link. In section \ref{sec_twoblock}, we construct the
two-block Gibbs sampler for the Bayesian probit linear mixed model
under improper priors. In section \ref{sec_convergence}, we prove
geometric ergodicity of the underlying Markov chain.  In section
\ref{sec_px-da}, we present a corresponding Haar PX-DA
algorithm. Section \ref{sec_conclusion} contains some conclusions and
discussions. Finally, the proofs of posterior propriety and geometric
convergence of the Gibbs sampler appear in the appendices.

\section{Propriety of posterior distributions}
\label{sec_propriety}

In this section, we discuss conditions under which the posterior
density \eqref{eq:postbetphi} is proper. The results in this section
hold for GLMMs with a general link function. Let $F(\cdot)$ be a cumulative distribution function, and
consider the link function $F^{-1}(\cdot)$. Thus instead of the probit linear mixed model in
\eqref{eq_link}, in this section we consider a GLMM with
\begin{equation}
\label{eq_link_general}
P(Y_i =1) = F(\bm{x}_i^T\bm{\beta} + \bm{z}_i^T\bm{u}).
\end{equation}

Posterior propriety for Bayesian GLMMs under improper priors has been
discussed in \cite{chen2002necessary}. We will first describe
\pcite{chen2002necessary} conditions. Then we will show, through
examples, that these conditions often do not hold in
practice. Finally, our conditions for posterior propriety will be
presented. 

Let $c_i = 1$ if $y_i = 0$ and $c_i = -1$ if $y_i = 1$ for
$i=1,\dots, n$. Suppose $\bm{W}_{n\times(p+q)}^*$ is a matrix whose
$i$th row is $c_i(\bm{x}_i^T,\bm{z}_i^T)$. In the special case when
$b_j=0$, that is, when $\tau_j$ has the power prior
$\pi(\tau_j) \propto \tau_j^{a_j-1}$ for $j=1,\dots,r$, a
straightforward extension of \pcite{chen2002necessary} Theorem 4.2
shows that the corresponding posterior distribution is proper if the
following conditions hold:
\begin{enumerate}[label={(A\arabic*)}]
\item  $\bm{W} = (\bm{X}, \bm{Z})$ is a full rank matrix;
\item There exists an $n\times 1$ positive vector $\bm{e}>0$ such that $\bm{e}^T\bm{W}^* = 0$;
\item $2a_j+ q_j >0$ for $j=1,2,\dots,r$;
\item $a_j < 0 $ for $j=1,\dots, r$;
\item $E \vert \delta \vert^{p - 2\sum_{j=1}^ra_j} <\infty$, where $\delta \sim F$.
\end{enumerate}
\cite{roy2007convergence} provided a simple method for checking the condition A2 using publicly available softwares. 

The condition A1 assumes that $\bm{W}$ is a full rank
matrix. Unfortunately, when $\bm{Z}$ is a design matrix with elements
1's and 0's, which is pretty common in practice, this assumption may
not hold. For example, we consider the following important generalized
two-way random effects model
\begin{equation}
\label{eq_link_example}
F^{-1}(P(Y_{ij} = 1)) = \beta + \alpha_i + \gamma_j,
\end{equation}
for $i = 1,2, \dots, n_1$, $j = 1,2,\dots, n_2$. Here, the $\alpha_i$'s
are i.i.d $N(0,1/\tau_1)$, and the $\gamma_j$'s are i.i.d
$N(0,1/\tau_2)$. There are total $n = n_1\times n_2$ observations and
we order them as
$\bm{Y} = (Y_{11},\dots, Y_{1n_2}, \dots, Y_{n_11}, \dots,
Y_{n_1n_2})$. In this example, $p=1$, and $\bm{X} = \bm{1}_n$ is an
$n\times 1$ column vector of ones. Also, there are $r=2$ random effects
with $q_1 = n_1$, $q_2 = n_2$, $q = q_1 + q_2 = n_1 + n_2$, and
$\bm{Z} = (\bm{Z}_1, \bm{Z}_2 )$, where
$\bm{Z}_1 = \bm{I}_{n_1} \otimes \bm{1}_{n_2}$ and
$\bm{Z}_2 = \bm{1}_{n_1}\otimes \bm{I}_{n_2}$ with $\otimes$ denoting the
Kronecker product.  It can be checked that the rank of
$\bm{W} = (\bm{X},\bm{Z})$ is $n_1 + n_2 - 1$. Thus $\bm{W}$ is not a
full rank matrix. 

We now provide Theorem \ref{them3} showing the posterior propriety
without the assumption A1. We also consider the more general prior
$\pi(\tau_j)$ given in \eqref{eq:improper_prior}, that is, $b_j$ may
not be zero. We use certain transformations of the regression
parameters $\bm{\beta}$ and random effects $\bm{u}$ to circumvent the
problem with non-full rank matrix $\bm{W}$. Assume that the first
column of $\bm{X}$ is a vector of 1's corresponding to an intercept
term $\beta_0$ in
$\bm{\beta} = (\beta_0, \beta_1,\dots, \beta_{p-1})^T$. Let
$\bm{Z} = (\bm{Z}_1,\dots,\bm{Z}_r)$, where $\bm{Z}_j$ is an
$n\times q_j$ matrix such that the $(ik)$th element is 1 if the
observation $i$ is observed at the $k$th level of the random effect
$\bm{u}_{j}=\left(u_{j1},\dots,u_{jq_{j}}\right)^{T}$, 0 otherwise,
for $i=1, \dots, n$, $k=1,\dots, q_j$ and $j=1,\dots, r$.  Consider
the following transformations,
\begin{eqnarray}
\label{eq_trans1}
\mu_{0} &=&\beta_{0}+\sum_{j=1}^{r}u_{j1},\\
\label{eq_trans2}
d_{jk} &=&u_{j,k+1}-u_{j1}, \text{for } k=1,\dots,q_{j}-1,j=1,\dots,r.
\end{eqnarray}
Thus $\mu_0$ is the sum of the intercept term and the first level
effect of all $r$ random effects. Also the (transformed) random
effects $d_{jk}$'s denote the differences of the random effect
compared to the first level effect.

Let
$\tilde{\bm{\eta}}=\left(\mu_{0},\beta_{1,},\dots,\beta_{p-1},d_{11},\dots,d_{1,q_{1}-1},\dots
  d_{r1},\dots,d_{r,q_{r}-1}\right)^{T}$.  Define
$\tilde{\bm{Z}} = (\tilde{\bm{Z}}_1,\dots, \tilde{\bm{Z}}_r)$, where
the $n\times (q_j -1)$ matrix $\tilde{\bm{Z}}_j$ is $\bm{Z}_j$ without
its first column. Thus, the vector $\bm{W}\bm{\eta}$ is the same as
the vector $\tilde{\bm{W}}\tilde{\bm{\eta}}$, where
$\tilde{\bm{W}}=\left(\bm{X},\tilde{\bm{Z}}\right)$ with $i$th row
$\tilde{\bm{w}}_i^T$. Let $\tilde{\bm{W}}^*$ be a matrix whose $i$th
row is
$ \tilde{\bm{w}}_i^{*T} =c_i\tilde{\bm{w}}_i^T = c_i(\bm{x}_i^T,
\tilde{\bm{z}}_i^T)$, where $\tilde{\bm{z}}_{i}^T$ is the $i$th row of
$\tilde{\bm{Z}}$. 

For the example \eqref{eq_link_example}, the
transformed parameters $\mu_0$ and $d_{jk}$'s become
\begin{eqnarray*}
\mu_0 & =& \beta + \alpha_1 + \gamma_1,\\
d_{1k} &=& \alpha_{k+1} - \alpha_1 \text{ for } k=1,\dots,n_1 - 1,\\
d_{2k} &=& \gamma_{k+1} - \gamma_1 \text{ for } k=1,\dots,n_2 - 1.
\end{eqnarray*}
Thus in this example, we have $\tilde{\bm{\eta}} = (\mu_0, d_{11},\dots, d_{1,n_1-1}, d_{21},\dots, d_{2,n_2-1})^T$. Also note that $\tilde{\bm{W}}$ is a full rank matrix in this example, although $\bm{W}$ is not.

\begin{theorem}
\label{them3}
Assume the following conditions hold,
\begin{enumerate}[label={(B\arabic*)}]
\item $a_j < b_j = 0$, $q_j \geq 2$ or $b_j >0$ for $j=1,\dots,r$;
\item $2a_j + q_j -1>0$ for $j=1,\dots,r$;
\item $\tilde{\bm{W}}$ is a full rank matrix;
\item There exists an $n\times 1$ positive vector $\bm{e}>0$ such that $\bm{e}^T\tilde{\bm{W}}^* = 0$.

\item $E \vert \delta \vert ^{p+t} <\infty$, where $t=\sum_{j=1}^r [- 2a_j I(b_j = 0) + (q_j -1) I(b_j>0)]$, and $\delta \sim F$.
%If $a_j < b_j = 0$, $E \vert \delta \vert ^{p - \sum_{j=1}^r a_j } < \infty$, and if $b_j >0$, $E \vert \delta \vert ^{p+q-r} <\infty$, where $\delta \sim F$;

Then the joint posterior density \eqref{eq:postbetphi} corresponding the GLMM \eqref{eq_link_general} is proper, i.e,
\begin{align}
\label{eq:integral}
 & \int_{\mathbb{R}_{+}^{r}}  \int_{\mathbb{R}^{q}} \int_{\mathbb{R}^{p}}\prod_{i=1}^{n}\left[ F\left(\bm{x}_{i}^{T}\bm{\beta}+\bm{z}_{i}^{T}\bm{u}\right)\right]^{y_{i}} \left[1-F\left(\bm{x}_{i}^{T}\bm{\beta}+\bm{z}_{i}^{T}\bm{u}\right)\right]^{1-y_{i}} \nonumber\\
 & \cdot\prod_{j=1}^{r}\tau_{j}^{\frac{q_{j}}{2}+a_j-1}\exp\left[-\tau_j \left( b_j + \frac{1}{2}\bm{u}_j^T\bm{u}_j\right)\right]d\bm{\beta}d\bm{u}d\bm{\tau}<\infty.
\end{align}
\end{enumerate}
\end{theorem}

A proof of Theorem \ref{them3} is given in Appendix \ref{proof_thm3}.

\begin{remark}
  When probit link is considered, that is $F(\cdot) = \Phi(\cdot)$,
  the moment condition $\textit{B5}$ holds automatically. Thus, for
  probit linear mixed models, the posterior density
  \eqref{eq:postbetphi} is proper under $\textit{B1} - \textit{B4}$.
\end{remark}

\section{A two-block Gibbs sampler}
\label{sec_twoblock}

We begin with deriving the joint density
$\pi(\bm{\beta}, \bm{u},\bm{v},\bm{\tau}|\bm{y})$ mentioned in the
introduction. Define the joint posterior density (up to a normalizing
constant) of $\bm{\beta}, \bm{u}, \bm{v}, \bm{\tau}$, if it exists, as
 \begin{align}
\label{eq_joint_improper}
\pi(\bm{\beta},\bm{u}, \bm{v}, \bm{\tau}|\bm{y}) & \propto  \prod_{i=1}^n \exp\left\{-\frac{1}{2}\left( v_i - \bm{x}_i^T\bm{\beta} - \bm{z}_i^T\bm{u} \right)^2 \right\} \nonumber\\ 
& \times  \prod_{i=1}^n  \left[1_{\text{(0,\ensuremath{\infty})}}\left(v_{i}\right)\right]^{y_{i}}\left[1_{\left(-\infty,0\right]}\left(v_{i}\right)\right]^{1-y_{i}} \nonumber \\
& \times \prod_{j=1}^r \tau_j^{\frac{q_j}{2}+a_j-1}\exp\left\{-\tau_j\left( b_j + \frac{\bm{u}_j^T\bm{u_j}}{2}\right)\right\}.
\end{align}
From \eqref{eq:lik} and \eqref{eq:vy} it follows that
\eqref{eq:margjt} holds. In section \ref{sec_propriety}, we discussed
conditions under which the posterior density
$\pi(\bm{\beta},\bm{\tau}\vert \bm{y})$ given in \eqref{eq:postbetphi}
and hence the joint posterior density \eqref{eq_joint_improper} is
proper. Note that, these posterior densities are proper if and only if
$c(\bm{y}) < \infty$.

Standard calculations show that the conditional density of $\bm{\eta}$ is
\begin{equation}
\label{eq_posterior1}
\pi(\bm{\eta}|\bm{v},\bm{\tau},\bm{y}) \propto  \exp\left[-\frac{1}{2}\left(\bm{v}-\bm{W}\bm{\eta}\right)^{T}\left(\bm{v}-\bm{W}\bm{\eta}\right)\right] \cdot\exp\left[-\frac{1}{2} \bm{u}^{T} \bm{D}(\bm{\tau}) \bm{u} \right],
\end{equation}
Thus,
\begin{equation}
\label{eq:eta}
\bm{\eta}|\bm{v}, \bm{\tau},\bm{y}\sim N_{p+q}\left(\bm{\bm{\Sigma}}^{-1}\bm{W}^T\bm{v},\bm{\Sigma}^{-1}\right),
\end{equation}
where
\begin{equation}
\label{eq:improper_eta_sigma}
\bm{\Sigma}=\left(\begin{array}{cc}
\bm{X}^{T}\bm{X} & \bm{X}^{T}\bm{Z}\\
\bm{Z}^{T}\bm{X} & \bm{Z}^{T}\bm{Z}+\bm{D}(\bm{\tau})
\end{array}\right). 
\end{equation}

Similarly, the conditional density of $(\bm{v}, \bm{\tau})$ is
\begin{eqnarray*}
  \pi(\bm{v},\bm{\tau}|\bm{\eta},\bm{y}) & \propto & \prod_{i=1}^{n} \phi \left(v_i- \bm{w}_i^T\bm{\eta};0,1\right)\left[1_{\text{(0,\ensuremath{\infty})}}\left(v_{i}\right)\right]^{y_{i}}\left[1_{\left(-\infty,0\right]}\left(v_{i}\right)\right]^{1-y_{i}}\\
                                         &  & \times\prod_{j=1}^{r}\tau_{j}^{\frac{q_{j}}{2} + a_j - 1}\exp\left[-\tau_{j} \left(b_j + \frac{1}{2}\bm{u}_{j}^{T}\bm{u}_{j}\right) \right] ,
\end{eqnarray*}
where $\bm{w}_i^T$ is the $i$th row of $\bm{W}$ for $i=1,\dots,n$. Thus, conditional on $(\bm{\eta}, \bm{y})$, $v_i$, $i = 1, \dots, n$ and $\bm{\tau}$ are independent. We have
\begin{equation}
\label{eq:v}
v_i| \bm{\eta}, \bm{y} \overset{\text{ind}}\sim \text{TN}(\bm{w}_i^T\bm{\eta},1,y_i), \, i=1,\dots, n,
\end{equation}
where $\text{TN}(\mu, \sigma^2, \omega)$ denotes the distribution of
the normal random variable with mean $\mu$ and variance $\sigma^2$,
that is truncated to have only positive values if $\omega = 1$, and
only nonpositive values if $\omega = 0$. Also conditional on
$\bm{\eta}, \bm{y}$, $\tau_j$'s are independent with
$\tau_j \sim \text{Gamma}\left(a_j+q_j/2,
  b_j+\bm{u}_j^T\bm{u}_j/2\right)$ for $j=1,\dots,r$.

Thus, one single iteration of the block Gibbs sampler $\{\bm{\eta}^{(m)}, \bm{v}^{(m)}, \bm{\tau}^{(m)}\}_{m=0}^{\infty}$ has the following two steps:

\begin{algorithm}[H]
\caption{The $(m+1)$st iteration of the two-block Gibbs sampler}
\label{algorithm2}
\begin{algorithmic}[1]
\STATE Draw $\tau_j^{(m+1)} \overset{\text{ind}}\sim \text{Gamma}\left(a_j+q_j/2, b_j+\bm{u}_j^{(m)T}\bm{u}_j^{(m)}/2\right)$ for $j=1,\dots, r$, and independently draw $v_i^{(m+1)} \overset{\text{ind}}\sim \text{TN}(\bm{w}_i^T\bm{\eta}^{(m)},1,y_i), \, i=1,\dots, n$.
%\STATE[1]
%\hspace{.18in}
\STATE Draw $\bm{\eta}^{(m+1)} \sim N_{p+q}\left(\left[\bm{\Sigma}^{(m+1)}\right]^{-1}\bm{W}^{T}\bm{v}^{(m+1)},\left[\bm{\Sigma}^{(m+1)}\right]^{-1}\right)$, where $\bm{\Sigma}^{(m+1)}$ is evaluated at $\bm{\tau}^{(m+1)}$.
\end{algorithmic}
\end{algorithm}

\section{Geometric ergodicity of the block Gibbs sampler}
\label{sec_convergence}
In this section, we establish the geometric rate of convergence of the block
Gibbs sampler
$\{\bm{\eta}^{(m)}, \bm{v}^{(m)}, \bm{\tau}^{(m)}\}_{m=0}^{\infty}$.
Since it is a two-block Gibbs sampler, it has the same rate of convergence
as the $\bm{\eta}$-marginal Markov chain $\{\bm{\eta}^{(m)}\}_{m=0}^{\infty}$
\citep{robe:rose:2001}. Below we analyze this $\bm{\Psi} \equiv \{\bm{\eta}^{(m)}\}_{m=0}^{\infty}$ chain.

 Let $\bm{\eta}^{\prime}$ be the current state and $\bm{\eta}$ be the
 next state of the Markov chain $\bm{\Psi}$, then the Markov transition density (Mtd) of $\bm{\Psi}$ is
\begin{equation}
\label{eq:mtd}
 k(\bm{\eta}|\bm{\eta}^{\prime}) =  \int_{\mathbb{R}_+^r} \int_{\mathbb{R}^n} \pi(\bm{\eta}|\bm{v},\bm{\tau},\bm{y})\pi(\bm{v},\bm{\tau}|\bm{\eta}^{\prime},\bm{y})d\bm{v}d\bm{\tau},
 \end{equation}
 where $\pi(\cdot|\cdot,\bm{y})$'s are the conditional densities from
 section~\ref{sec_twoblock}. Routine calculations show that
 $k(\bm{\eta}|\bm{\eta}^{\prime})$ is reversible and thus is invariant
 with respect to the marginal density of $\bm{\eta}$ denoted as
 $\pi(\bm{\eta}|\bm{y}) \equiv \int_{\mathbb{R}_+^r}
 \int_{\mathbb{R}^{n}}
 \pi(\bm{\eta},\bm{v},\bm{\tau}|\bm{y})d\bm{v}d\bm{\tau}$. Let
 $h: \mathbb{R}^{p+q} \mapsto \mathbb{R}$ be a real valued
 function. Suppose our interest is to estimate the (posterior) mean
 $ E (h(\bm{\eta})|\bm{y}) \equiv \int_{\mathbb{R}^{p+q}}
 h(\bm{\eta})\pi(\bm{\eta}|\bm{y}) d\bm{\eta}$. Since
 $k(\bm{\eta}|\bm{\eta}^{\prime})$ is strictly positive, the Markov
 chain $\bm{\Psi}$ is Harris ergodic \citep{meyn1993markov}. Thus if
 $E (|h(\bm{\eta})| |\bm{y}) < \infty$, then $E(h(\bm{\eta})|\bm{y})$
 can be consistently estimated by
\begin{equation*}
\bar{h}_m = \frac{1}{m}\sum_{i=0}^{m-1}h(\bm{\eta}^{(i)}).
\end{equation*}
As mentioned in the introduction, in order to provide an
asymptotically valid confidence interval for $E(h(\bm{\eta})|\bm{y})$
based on $\bar{h}_m$, we need to establish a CLT for $\bar{h}_m$. We
say a CLT exists for $\bar{h}_m$ if there exists a constant
$\sigma^2_h \in (0,\infty)$ such that,
\begin{equation}
\label{eq:clt}
\sqrt{m}\left(\bar{h}_m - E(h(\bm{\eta})|\bm{y})\right)\overset{d}\rightarrow N\left(0,\sigma_h^2\right) \, \text{as } m\rightarrow \infty .
\end{equation}
If \eqref{eq:clt} holds, and a consistent estimator $\hat{\sigma}_h^2$
of $\sigma_h^2$ is available, then the standard errors
$\hat{\sigma}_h/\sqrt{m}$ can be used to provide an asymptotic
confidence interval for $E(h(\bm{\eta})|\bm{y})$
\citep{roy2007convergence}. Unfortunately, Harris ergodicity of
$\bm{\Psi}$ does not guarantee \eqref{eq:clt}, although it ensures
consistency of $\bar{h}_m$. One method of proving \eqref{eq:clt} is to
establish the geometric rate of convergence for the Markov chain
$\bm{\Psi}$ \citep{jones2001honest}.  Geometric ergodicity of $\bm{\Psi}$ also
allows for consistent estimation of $\sigma_h^2$ using batch means or
spectral variance methods \citep{flegal2010batch}.

Let $\mathscr{B}$ denote the Borel $\sigma$-algebra of $\mathbb{R}^{p+q}$ and
$K(\cdot,\cdot)$ be the Markov transition function corresponding to the Mtd
$k(\cdot, \cdot)$ in \eqref{eq:mtd}, that is, for any set
$O \in \mathscr{B}$, $\bm{\eta}^{\prime} \in \mathbb{R}^{p+q}$ and any
$j=0,1,\dots,$
\begin{equation}
K(\bm{\eta}^{\prime}, O) = \mbox{Pr}(\bm{\eta}^{(j+1)} \in O| \bm{\eta}^{(j)} = \bm{\eta}^{\prime}) = \int_{O} k(\bm{\eta}| \bm{\eta}^\prime) d\bm{\eta}.
\end{equation}
Then the $m$-step Markov transition function is
$K^m(\bm{\eta}^{\prime}, O) = \mbox{Pr}(\bm{\eta}^{(m+j)} \in O |
\bm{\eta}^{(j)} = \bm{\eta}^{\prime})$. Let $\Pi(\cdot|\bm{y})$ be the
probability measure with density $\pi(\bm{\eta}|\bm{y})$. The Markov
chain $\bm{\Psi}$ is geometrically ergodic if there exists a constant
$0 < t <1$ and a function $J: \mathbb{R}^{p+q} \mapsto \mathbb{R}^+$
such that for any $\bm{\eta} \in \mathbb{R}^{p+q}$,
\begin{equation}
\label{eq:ge}
||K^m(\bm{\eta},\cdot) - \Pi(\cdot|\bm{y})||_{\text{TV}}:=\sup_{O\in \mathscr{B}} |K^m(\bm{\eta}, O) - \Pi(O|\bm{y})| \leq J(\bm{\eta}) t^m.
\end{equation}
Harris ergodicity of $\bm{\Psi}$ implies that
$||K^m(\bm{\eta},\cdot) - \Pi(\cdot|\bm{y})||_{\text{TV}} \downarrow 0$ as
$m \rightarrow \infty$, while \eqref{eq:ge} guarantees its exponential
rate of convergence. \cite{roberts1997geometric} showed that since
$\bm{\Psi}$ is reversible, if \eqref{eq:ge} holds then there exists a CLT,
that is \eqref{eq:clt} holds, for all $h$ with
$E (h^2(\bm{\eta})|\bm{y}) < \infty$. 

In section \ref{sec_propriety}, we provided two sets of conditions for
posterior propriety. While the first set of conditions
$(\textit{A1} - \textit{A5})$ holds in the special case $b_j =0$ for
all $j=1,\dots,r$, Theorem \ref{them3} holds for the general prior
$\pi(\tau_j)$ given in \eqref{eq:improper_prior}. In Theorems
\ref{them2} and \ref{them4}, we provide conditions under which the
Markov chain $\bm{\Psi}$ is geometrically ergodic, that is,
\eqref{eq:ge} holds. Here we consider the general form of the prior
distribution of $\tau_j$ as given in \eqref{eq:improper_prior}. Thus
the parameters $b_j$'s are not assumed to be zero.  Since geometric
ergodicity implies posterior propriety, Theorem~\ref{them2} also
provides conditions for posterior propriety for the probit linear
mixed models in the general case when $b_j \neq 0$.

\begin{theorem}
\label{them2}
The Markov chain underlying the block Gibbs sampler is geometrically ergodic if the following conditions hold:
\begin{enumerate}[label={(\arabic*)}]
\item $a_j <b_j  = 0$ or $b_j>0$ for $j=1,\dots,r$;
\item $(A1) - (A3)$ hold.
\end{enumerate}
\end{theorem}

A proof of Theorem \ref{them2} is given in the Appendix
\ref{proof_thm2}. Theorem \ref{them4} shows geometric convergence of
the Markov chain underlying the Gibbs sampler given in Algorithm
\ref{algorithm2} without the assumption A1.
\begin{theorem}
\label{them4}
The block Gibbs sampler is geometrically ergodic under the following conditions:
\begin{enumerate}[label={(\arabic*)}]
\item $(B1) - (B4)$ hold;
\item 
There exists an $s \in (0,1] \cap (0,\tilde{s})$  such that 
\begin{equation}
\label{eq_condition_thm4}
2^{-s}\sum_{j=1}^{r}\frac{\Gamma\left(q_{j}/2+a_{j}-s\right)}{\Gamma\left(q_{j}/2+a_{j}\right)}\left[tr\left(R_{j}\left(\bm{I}-P_{\bm{Z}^{T}\left(\bm{I}-P_{X}\right)\bm{Z}}\right)R_{j}^{T}\right)\right]^{s}<1,
\end{equation}
where $\tilde{s} = \min\{a_1 + q_1/2,\dots, a_r + q_r/2 \}$, $R_{j}$ is a $q_{j}\times q$ matrix with 0's and 1's such that $R_{j}\bm{u}=\bm{u}_{j}$ and $P_{\bm{Z}^{T}\left(\bm{I}-P_{X}\right)\bm{Z}}$ is the projection matrix on the column space of $\bm{Z}^{T}\left(\bm{I}-P_{X}\right)\bm{Z}$.
\end{enumerate}
\end{theorem}

A proof of Theorem \ref{them4} is given in Appendix \ref{proof_thm4}.

\begin{remark}
  The extra condition (2) in Theorem $\ref{them4}$ compared to Theorem
  $\ref{them2}$ is due to the lack of the full rank assumption of
  $\bm{W}$, and the need to include an extra term in the drift
  function used to prove Theorem \ref{them4}. This condition is also
  used in \cite{roman2012convergence}, who provide some discussions on
  this. The left-hand side of \eqref{eq_condition_thm4} can be
  evaluated at values of $s$ on a fine grid in the interval
  $ (0,1] \cap (0,\tilde{s})$ to numerically check the condition. Note
  that, $R_j$ is the matrix that extracts $\bm{u}_j$ out of
  $\bm{u}$. Thus when $r >1$,
  $tr\left(R_{j}\left(\bm{I}-P_{\bm{Z}^{T}\left(\bm{I}-P_{X}\right)\bm{Z}}\right)R_{j}^{T}\right)$
  is the sum of the $q_j$ diagonal elements of
  $\bm{I}-P_{\bm{Z}^{T}\left(\bm{I}-P_{X}\right)\bm{Z}}$ corresponding
  to the $j$th random effect.
\end{remark}

\section{A Haar PX-DA algorithm}
\label{sec_px-da}
As mentioned in section \ref{sec_introduction}, DA algorithms often suffer from slow
convergence and high autocorrelations. \cite{liu1999parameter}
proposed parameter expansion for data augmentation (PX-DA) algorithms
for speeding up the convergence of DA
algorithms. \cite{hobert2008theoretical} compared the performance of
PX-DA algorithms based on a Haar measure (called Haar PX-DA algorithms)
with PX-DA algorithms based on a probability measure and DA
algorithms. In particular, they showed that, under some mild
conditions, the Haar PX-DA algorithms are better than the general PX-DA
algorithms and the DA algorithms in both the efficiency ordering and the operator norm
ordering. As shown in \cite{hobert2008theoretical}, compared to the DA
algorithm, in PX-DA, an extra step is added (sandwiched) between the two steps of
the original DA algorithm.
In order to construct this extra step, we
derive the marginal density
\begin{eqnarray}
\label{eq:margvtau1}
  \pi\left(\bm{v}, \bm{\tau}|\bm{y}\right) & = & \int_{\mathbb{R}^{p+q}}\pi(\bm{\eta},\bm{v},\bm{\tau}|\bm{y})d\bm{\eta}\\
                                           & \propto & \prod_{i=1}^{n}\left[1_{\text{(0,\ensuremath{\infty})}}\left(v_{i}\right)\right]^{y_{i}}\left[1_{\left(-\infty,0\right]}\left(v_{i}\right)\right]^{1-y_{i}} \prod_{j=1}^{r}\tau_{j}^{\frac{q_{j}}{2}+a_j} e^{-b_j\tau_j}\nonumber\\
                                           & & \cdot \vert \bm{\Sigma}\vert^{-1/2}\exp\left\{ -\frac{1}{2}\bm{v}^{T}M_1\bm{v}\right\} , \nonumber
\end{eqnarray}
where  $M_1 = \bm{I}-\bm{W}\bm{\Sigma}^{-1}\bm{W}^{T}$. 

Let $\mathcal{Z}$ denote the subset of $\mathbb{R}^n$ where $\bm{v}$
lives, that is, $\mathcal{Z}$ is the Cartesian product of $n$ half
(positive or nonpositive) lines, where the $i$th component is
$(0, \infty)$ (if $y_i =1$) or $(-\infty, 0]$ (if $y_i =0$).  Let $G$
be the unimodular multiplicative group on $\mathbb{R}_+$ with Haar
measure $\nu(dg) = dg/g$, where $dg$ is the Lebesgue measure on
$\mathbb{R}_+$.  For constructing an efficient extra step, as in
\cite{roy2014efficient}, we let the group $G$ act on
$\mathcal{Z} \times \mathbb{R}_+^r$ through a group action
$T(\bm{v}, \bm{\tau}) = (g\bm{v},\bm{\tau}) = (gv_1, gv_2, \dots,
gv_n, \bm{\tau})$. With the group action defined this way, it can be
shown that the Lebesgue measure on $\mathcal{Z} \times \mathbb{R}_+^r$
is relatively left invariant with multiplier $\chi (g) = g^n$
\citep{roy2014efficient, hobert2008theoretical}. Following
\cite{hobert2008theoretical}, consider a probability density function
$\vartheta(g)$ on $G$ where
\begin{equation}
\vartheta\left(g\right) dg \propto \pi\left(g\bm{v}, \bm{\tau}|\bm{y}\right) \chi \left( g \right) \nu(dg) \propto  g^{n-1}\exp\left\{ -\frac{1}{2}g^{2}\bm{v}^{T} M_1 \bm{v} \right\} dg.
\label{eq_g1}
\end{equation}
Since propriety of the posterior density \eqref{eq_joint_improper}
implies that $\pi(\bm{v}, \bm{\tau}|\bm{y})$ is a valid density,
$\bm{v}^TM_1\bm{v}$ can be zero only on a set of measure 0 (in
$\bm{v}$). Thus given $(\bm{v},\bm{\tau})$, $\vartheta\left(g\right)$
is a valid density. From \cite{hobert2008theoretical}, it follows that
the transition
$(\bm{v}, \bm{\tau}) \rightarrow (\bm{v}', \bm{\tau}) \equiv T(\bm{v},
\bm{\tau}) = (g\bm{v},\bm{\tau})$ where $g \sim \vartheta(g)$, is
reversible with respect to $\pi(\bm{v},\bm{\tau}|\bm{y})$ defined in
\eqref{eq:margvtau1}. Given $\bm{\eta}^{(m)}$, below are the three
steps involved in the $(m+1)$st iteration of the Haar PX-DA algorithm
to move to the new state $\bm{\eta}^{(m+1)}$.

\begin{algorithm}[H]
\caption{The $(m+1)$st iteration of the Haar PX-DA algorithm}
\label{algorithm_g1}
\begin{algorithmic}[1]
  \STATE $\tau_j \sim \text{Gamma}\left(a_j+q_j/2, b_j+\bm{u}_j^{(m)T}\bm{u}_j^{(m)}/2\right)$, for $j=1,\dots, r$ and independently draw $v_i| \bm{\eta}^{(m)}, \bm{y} \overset{\text{ind}}\sim$ $
  \text{TN}(\bm{w}_i^T\bm{\eta}^{(m)},1,y_i)$ for  $i=1,\dots, n$.

\STATE Draw $g^2$ from $\text{Gamma}(n/2, \bm{v}^TM_1\bm{v}/2)$.

\STATE Set $v_i^{\prime} = gv_i$ and let $\bm{v}^{\prime} = (v_1^{\prime},\dots, v_n^{\prime})^T$. Draw 
\[\bm{\eta}^{(m+1)} \sim N_{p+q}\left(\bm{\Sigma}\left(\bm{\tau}\right)^{-1}\bm{W}^{T}\bm{v}^{\prime },\bm{\Sigma}\left(\bm{\tau}\right)^{-1}\right).\]
 \end{algorithmic}
\end{algorithm}
The Mtd of the above Haar PX-DA algorithm can be written as
\begin{equation}
\label{eq:mtdpxda}
  k^*\left(\bm{\eta}|\bm{\eta}^{\prime}\right) = \int_{\mathbb{R}^n}  \int_{\mathbb{R}_{+}^{r}} \int_{\mathbb{R}^n}\pi\left(\bm{\eta}|\bm{v}^{\prime},\bm{\tau},\bm{y}\right) Q\left(\bm{v},d\bm{v}^{\prime}\right) \pi\left(\bm{v},\bm{\tau}|\bm{\eta}^{\prime},\bm{y}\right)d\bm{v}d\bm{\tau} d\bm{v}^{\prime},
\end{equation}
where $Q(\cdot,\cdot)$ is the Markov transition function corresponding
to the move
$(\bm{v}, \bm{\tau}) \rightarrow (\bm{v}', \bm{\tau}) = T(\bm{v},
\bm{\tau})$.  Let $K^*$ and $K$ be the Markov operators associated
with the Mtds $k^*$ and $k$ defined in \eqref{eq:mtdpxda} and
\eqref{eq:mtd} respectively. From \cite{hobert2008theoretical}, we
have $\Vert K^*\Vert_{\text{OP}}\leq \Vert K\Vert_{\text{OP}}$, where $\Vert K\Vert_{\text{OP}}$ denotes the norm of the operator $K$ \cite[see also][]{roy:2012a}. Since the block Gibbs sampler is
geometrically ergodic, we have $\Vert K^*\Vert_{\text{OP}}\leq \Vert K\Vert_{\text{OP}} < 1$ \citep{roberts1997geometric}.  Thus we
have the following corollary.

\begin{corollary} Under the conditions of Theorem $\ref{them2}$ or Theorem $\ref{them4}$, the Markov chain  underlying the Haar PX-DA algorithm described in Algorithm $\ref{algorithm_g1}$ is geometrically ergodic.
\end{corollary}

The extra step of Algorithm \ref{algorithm_g1} is a single draw
from the univariate density $\vartheta(g)$, which is easy to sample
from. Thus, the computational burden, per iteration, for the Haar
PX-DA algorithm is similar to that of the block Gibbs sampler
described in section~\ref{sec_twoblock}. Two other Haar
PX-DA algorithms can be constructed by using group actions
$T_1(\bm{v}, \bm{\tau}) = (\bm{v}, g\bm{\tau})$ and
$T_2(\bm{v}, \bm{\tau}) = (g\bm{v}, g\bm{\tau})$. However, the
corresponding $\vartheta(g)$'s are not easy to sample from, thus we do not consider them here.

\section{Discussion}
\label{sec_conclusion}

We develop a two-block Gibbs sampler for the Bayesian probit linear
mixed models under improper priors. The block Gibbs algorithm samples the fixed effects and the random effects jointly. We prove the geometric ergodicity of the two-block Gibbs sampler, which guarantees the existence of central limit theorems for MCMC estimators under a finite second moment condition. We also propose the corresponding Haar PX-DA algorithm. The Haar PX-DA algorithm not only improve the efficiency of the Gibbs sampler, but also
inherit their geometric convergence properties.

Another popular link function is the logit link
function. \cite{polson2013bayesian} proposed a DA algorithm for the
logistic regression model. \cite{choi2013} proved the uniform
ergodicity of this DA algorithm. As mentioned in
\cite{polson2013bayesian}, their DA algorithm can be extended to the
logistic linear mixed model. However, the convergence properties of
the corresponding Markov chain have not been studied, and can be a
topic for future research. Another future project can be deriving
similar extensions of the results in \cite{roy:2012b} for proving
geometric convergence of Gibbs samplers for robit linear mixed
models. 

\section*{Acknowledgment}
The authors thank two anonymous reviewers and an anonymous associate editor for several helpful comments and suggestions that led to a substantially improved revision of the paper.

\appendix
\renewcommand{\thesubsection}{ \Alph{subsection}}
\section*{Appendices}
\label{sec_appen}

\subsection{Proof of Theorem \ref{them3}}
\label{proof_thm3}
\begin{proof}
Using the transformation $(\bm{\beta}^T,\bm{u}^T)^T \rightarrow (u_{11},\dots, u_{r1},\tilde{\bm{\eta}}^T)^T$, the integral in \eqref{eq:integral} can be written as, 
\begin{align}
\label{eq_int}
& \int_{\mathbb{R}^{p+q-r}}\int_{\mathbb{R}_{+}^{r}}\int_{\mathbb{R}^{r}}\prod_{i=1}^{n}\left[F\left(\tilde{\bm{w}}_{i}^{T}\tilde{\bm{\eta}}\right)\right]^{y_{i}}\left[1-F\left(\tilde{\bm{w}}_{i}^{T}\tilde{\bm{\eta}}\right)\right]^{1-y_{i}} \\
& \prod_{j=1}^{r}\tau_{j}^{\frac{q_{j}}{2}+a_{j}-1}\exp\left[-\frac{\tau_{j}}{2}\left(u_{j1}^{2}+\sum_{k=1}^{q_{j}-1}\left(d_{jk}+u_{j1} \right)^{2} + 2b_j\right)\right]du_{11}\cdots du_{r1}d\bm{\tau}d\tilde{\bm{\eta}},\nonumber 
\end{align}
where $\tilde{\bm{w}}_i$ is defined in section \ref{sec_propriety}. Let $\bar{d}_j = \sum_{k=1}^{q_j - 1}d_{jk}$. Then \eqref{eq_int} becomes,

\begin{align} 
 & \int_{\mathbb{R}^{p+q-r}}\int_{\mathbb{R}_{+}^{r}}\prod_{i=1}^{n}\left[F\left(\tilde{\bm{w}}_{i}^{T}\tilde{\bm{\eta}}\right)\right]^{y_{i}}\left[1-F\left(\tilde{\bm{w}}_{i}^{T}\tilde{\bm{\eta}}\right)\right]^{1-y_{i}} \left(2\pi\right)^{\frac{r}{2}} \nonumber \\
 & \cdot\prod_{j=1}^{r}q_{j}^{-1/2}\tau_{j}^{\frac{q_{j}}{2}+a_{j}-\frac{3}{2}}\exp\left[-\frac{\tau_{j}}{2}\left(\sum_{k=1}^{q_{j}-1}\left(d_{jk}-\bar{d}_{j}\right)^{2}+\frac{q_{j}-1}{q_{j}}\bar{d}_{j}^{2} +2 b_j\right)\right]d\bm{\tau}d\tilde{\bm{\eta}}\nonumber \\
= & \int_{\mathbb{R}^{p+q-r}}\prod_{i=1}^{n}\left[F\left(\tilde{\bm{w}}_{i}^{T}\tilde{\bm{\eta}}\right)\right]^{y_{i}}\left[1-F\left(\tilde{\bm{w}}_{i}^{T}\tilde{\bm{\eta}}\right)\right]^{1-y_{i}} \left(2\pi\right)^{\frac{r}{2}}\nonumber \\
 & \cdot\prod_{j=1}^{r}q_{j}^{-1/2}\Gamma\left(q_{j}/2+a_{j}-1/2\right)2^{\frac{q_{j}}{2}+a_{j}-\frac{1}{2}} \nonumber \\
 &\cdot \left(\sum_{k=1}^{q_{j}-1}\left(d_{jk}-\bar{d}_{j}\right)^{2}+\frac{q_{j}-1}{q_{j}}\bar{d}_{j}^{2} + 2b_j\right)^{-\frac{q_{j}}{2}-a_{j}+\frac{1}{2}}d\tilde{\bm{\eta}}\nonumber \\
\leq & \varphi_{1}\int_{\mathbb{R}^{p+q-r}}\prod_{i=1}^{n}\left[F\left(\tilde{\bm{w}}_{i}^{T}\tilde{\bm{\eta}}\right)\right]^{y_{i}}\left[1-F\left(\tilde{\bm{w}}_{i}^{T}\tilde{\bm{\eta}}\right)\right]^{1-y_{i}} \nonumber \\
&\cdot\prod_{j=1}^{r}\left(\sum_{k=1}^{q_{j}-1}\left(d_{jk}-\bar{d}_{j}\right)^{2}+\frac{q_{j}-1}{q_{j}}\bar{d}_{j}^{2} + 2b_j\right)^{-\frac{q_{j}}{2}-a_{j}+\frac{1}{2}}d\tilde{\bm{\eta}} \label{eq:improper_integral},
\end{align}
where $\varphi_1$ is a constant depending on  $r$, $q_j$ and $a_j$, $j=1,\dots, r$.

Let $\delta_{i},i=1,\dots,n$ be $n$ i.i.d random variables with
distribution function $F$. Let
$\bm{\delta}^{*}=\left(c_{1}\delta_{1},\dots,c_{n}\delta_{n}\right)^{T}$,
where $c_i = 1 - 2y_i$ as defined in section \ref{sec_propriety}. We
have
$E\left[1\left\{ c_{i}\tilde{\bm{w}}_{i}^{T}\tilde{\bm{\eta}}\leq
    c_{i}\delta_{i}\right\} \right]=\left[
  F\left(\tilde{\bm{w}}_{i}^{T}\tilde{\bm{\eta}}\right)\right]^{y_{i}}\left[1-F\left(\tilde{\bm{w}}_{i}^{T}\tilde{\bm{\eta}}\right)\right]^{1-y_{i}}$,
for $i = 1,\dots, n$. Thus
\begin{equation}
\label{eq:improper_integral_ineq}
\prod_{i=1}^{n}\left[F\left(\tilde{\bm{w}}_{i}^{T}\tilde{\bm{\eta}}\right)\right]^{y_{i}}\left[1-F\left(\tilde{\bm{w}}_{i}^{T}\tilde{\bm{\eta}}\right)\right]^{1-y_{i}}=E\left[1\left\{ \tilde{\bm{W}}^{*}\tilde{\bm{\eta}}\leq\bm{\delta}^{*}\right\} \right],
\end{equation}
where $\tilde{\bm{W}}^{*}$ is the $n\times (p+q)$ matrix whose $i$th row is $c_i\tilde{\bm{w}}_{i}^{T}$.

Since conditions \textit{B3} and \textit{B4} are in force, according
to \cite{chen2001propriety} (Lemma 4.1), there exists a constant
$\varphi_0$ depending on $\tilde{\bm{W}}$ and $\bm{y}$, such that
$1\left\{
  \tilde{\bm{W}}^{*}\tilde{\bm{\eta}}\leq\bm{\delta}^{*}\right\}
\leq1\left\{ \left\Vert \tilde{\bm{\eta}}\right\Vert \leq
  \varphi_0\left\Vert \bm{\delta}^{*}\right\Vert \right\}$. Recall
that
$\tilde{\bm{\eta}} = (\mu_0, \beta_1,\dots, \beta_{p-1},
d_{11},\dots,\allowbreak d_{1,q_1-1},\dots, d_{r1}, \dots,
d_{r,q_r-1})^T = (\mu_0, \beta_1,\dots, \beta_{p-1}, \bm{d}_1^T,\dots,
\bm{d}_r^T)^T$, where
$\bm{d}_j = (d_{j1},\allowbreak \dots, d_{j, q_j - 1})^T$ for
$j=1,\dots, r$.  Thus from \eqref{eq:improper_integral} and
\eqref{eq:improper_integral_ineq} it follows that \eqref{eq_int} is
bounded above by
\begin{align}
\label{eq_int1}
 & \varphi_{1}E\left[\int_{\mathbb{R}^{p+q-r}}1\left\{ \left\Vert \tilde{\bm{\eta}}\right\Vert \leq \varphi_0\left\Vert \bm{\delta}^{*}\right\Vert \right\} \prod_{j=1}^{r}\left(\sum_{k=1}^{q_{j}-1}\left(d_{jk}-\bar{d}_{j}\right)^{2}+\frac{q_{j}-1}{q_{j}}\bar{d}_{j}^{2} + 2b_j\right)^{-\frac{q_{j}}{2}-a_{j}+\frac{1}{2}}d\tilde{\bm{\eta}}\right]  \nonumber \\
 & \leq2^{p}\varphi_0^{p}\varphi_{1}E\left[\left\Vert \bm{\delta}^{*}\right\Vert ^{p}\int_{A_d}\prod_{j=1}^{r}\left(\sum_{k=1}^{q_{j}-1}\left(d_{jk}-\bar{d}_{j}\right)^{2}+\frac{q_{j}-1}{q_{j}}\bar{d}_{j}^{2} + 2b_j\right)^{-\frac{q_{j}}{2}-a_{j}+\frac{1}{2}}d\bm{d}_{1}\cdots,d\bm{d}_{r}\right]\nonumber \\
 & \leq2^{p}\varphi_0^{p}\varphi_{1}E\left[\left\Vert \bm{\delta}^{*}\right\Vert ^{p}\int_{A_d}\prod_{j=1}^{r}\left(\frac{q_{j}-1}{q_{j}}\bar{d}_{j}^{2} + 2b_j\right)^{-\frac{q_{j}}{2}-a_{j}+\frac{1}{2}}d\bm{d}_{1}\cdots,d\bm{d}_{r}\right],
\end{align}
where $A_{d}=\left\{ \vert d_{jk}\vert\leq\varphi_{0}\Vert\bm{\delta}^{*}\Vert,j=1,\dots,r,\,k=1,\dots q_{j}-1\right\} $.

We consider two cases of condition \textit{B1} separately.

\textbf{Case 1}: $a_j <  b_j = 0 $, $q_j \geq 2$. 
If $q_j = 2$, we have 

\begin{align*}
\int_{\vert d_{j1}\vert\leq1}\left[\left(d_{j1}\right)^{2}\right]^{-a_{j}-\frac{1}{2}}dd_{j1} & =\left.-\frac{1}{2a_{j}}\left[\left(d_{j1}\right)^{2}\right]^{-a_{j}-\frac{1}{2}}d_{j1}\right|_{-1}^{1}\\
 & =-\frac{1}{2a_{j}}\left(1+1\right)=-\frac{1}{a_{j}}<\infty.
\end{align*}

For $q_j > 2$, note that,
\begin{align}
\label{eq_in2}
 & \int_{\left|d_{j1}\right|\leq1}\left[\left(\sum_{k=1}^{q_{j}-1}d_{jk}\right)^{2}\right]^{-\frac{q_{j}}{2}-a_{j}+\frac{1}{2}}dd_{j1} \nonumber \\
 &=\frac{1}{2-q_{j}-2a_{j}}\left.\left[\left(\sum_{k=1}^{q_{j}-1}d_{jk}\right)^{2}\right]^{-\frac{q_{j}}{2}-a_{j}+\frac{1}{2}}\left(\sum_{k=1}^{q_{j}-1}d_{jk}\right)\right|_{-1}^{1}  \nonumber \\
 & =\frac{1}{2-q_{j}-2a_{j}}\left[\left(1+\sum_{k=2}^{q_{j}-1}d_{jk}\right)^{2}\right]^{-\frac{q_{j}}{2}-a_{j}+\frac{1}{2}}\left(1+\sum_{k=2}^{q_{j}-1}d_{jk}\right) \nonumber\\
 & -\frac{1}{2-q_{j}-2a_{j}}\left[\left(-1+\sum_{k=2}^{q_{j}-1}d_{jk}\right)^{2}\right]^{-\frac{q_{j}}{2}-a_{j}+\frac{1}{2}}\left(-1+\sum_{k=2}^{q_{j}-1}d_{jk}\right) \nonumber\\
 & \leq\frac{1}{2-q_{j}-2a_{j}}\left\{ \left[\left(1+\sum_{k=2}^{q_{j}-1}d_{jk}\right)^{2}\right]^{-\frac{q_{j}}{2}-a_{j}+1}+\left[\left(-1+\sum_{k=2}^{q_{j}-1}d_{jk}\right)^{2}\right]^{-\frac{q_{j}}{2}-a_{j}+1}\right\}. 
\end{align}
Ignoring the constant multiple, continuing integrating \eqref{eq_in2}
with respect to $d_{j2},\dots, d_{j,q_j-1}$ consecutively, we arrive
at some linear combinations of terms
\begin{equation}
\label{eq_in3}
\left.\left[\left(\alpha_{0}+d_{j,q_{j}-1}\right)^{2}\right]^{-\frac{q_{j}}{2}-a_{j}+\frac{q_{j}-1}{2}}\left(\alpha_{0}+d_{j,q_{j}-1}\right)\right|_{-1}^{1},
\end{equation}
where $\alpha_0$'s are constants. Since $a_j < 0$, each of these terms in \eqref{eq_in3} is finite. Then
\begin{align}
 & \int_{A_{dj} }\left(\frac{q_{j}-1}{q_{j}}\bar{d}_{j}^{2}+b_{j}\right)^{-\frac{q_{j}}{2}-a_{j}+\frac{1}{2}}d\bm{d}_{j} \nonumber \\
= & \left[q_{j}\left(q_{j}-1\right)\right]^{\frac{q_{j}}{2}+a_{j}-\frac{1}{2}}\left(\varphi_0\left\Vert \bm{\delta}^{*}\right\Vert \right)^{-2a_{j}} \nonumber \\
 &\cdot\int_{\{\left|d_{jk}\right|\leq1, k = 1,\dots, q_j - 1\}}\left[\left(\sum_{k=1}^{q_{j}-1}d_{jk}\right)^{2}\right]^{-\frac{q_{j}}{2}-a_{j}+\frac{1}{2}}d\bm{d}_{j} \nonumber\\
\leq & \varphi_{2j}\left\Vert \bm{\delta}^{*}\right\Vert ^{-2a_{j}} \label{eq:improper_integral_c1},
\end{align}
where $A_{dj}=\left\{ \vert d_{jk}\vert\leq\varphi_{0}\Vert\bm{\delta}^{*}\Vert,k=1,\dots q_{j}-1\right\} $ and $\varphi_{2j}$ is a finite positive constant. 

\textbf{Case 2}: $  b_j > 0 $. 

We have
\begin{align}
 & \int_{A_{dj} }\left(\frac{q_{j}-1}{q_{j}}\bar{d}_{j}^{2}+2b_{j}\right)^{-\frac{q_{j}}{2}-a_{j}+\frac{1}{2}}d\bm{d}_{j}\leq  \int_{A_{dj}}(2b_{j})^{-\frac{q_{j}}{2}-a_{j}+\frac{1}{2}}d\bm{d}_{j} \nonumber \\
\leq & (2b_{j})^{-\frac{q_{j}}{2}-a_{j}+\frac{1}{2}}2^{q_{j}-1}\varphi_0^{q_{j}-1}\left\Vert \bm{\delta}^{*}\right\Vert ^{q_{j}-1}
\leq  \varphi_{3j}\left\Vert \bm{\delta}^{*}\right\Vert^{q_j-1}, \label{eq:improper_integral_c2}
\end{align}
where $\varphi_{3j}$ is a finite positive constant. 

Using \eqref{eq:improper_integral_c1}, \eqref{eq:improper_integral_c2},
and condition \textit{B5}, it follows that \eqref{eq_int} can be bounded above
by
\begin{align*}
 & 2^{p}\varphi_0^{p}\varphi_{1}E\left[\left\Vert \bm{\delta}^{*}\right\Vert ^{p}\prod_{j=1}^{r}\left\{ \varphi_{2j}\left\Vert \bm{\delta}^{*}\right\Vert ^{-2a_{j}}I(b_{j} = 0)+\varphi_{3j}\left\Vert \bm{\delta}^{*}\right\Vert ^{q_{j}-1}I(b_j >0)\right\} \right]\\
\leq & 2^{p}\varphi_0^{p}\varphi_{1}\prod_{j: b_j =0}\varphi_{2j} \prod_{j: b_j >0}\varphi_{3j} \left[E\left\Vert \bm{\delta}^{*}\right\Vert ^{p+\sum_{j=1}^{r}[-2a_{j}I(b_j = 0)+(q_j-1)I(b_j >0)]} \right]\\
 < &\infty.
\end{align*}
\end{proof}

\begin{remark}
  If $q_j = 1$ and $b_j = 0$, \eqref{eq_int} is $\infty$ since
  $\int_{\mathbb{R}_{+}}\tau_j^{q_j/2+a_j -3/2}d\tau_j = \infty$. If
  $b_j >0$, the posterior density \eqref{eq:postbetphi} can be
  proper even when $q_j=1$.
\end{remark}

\subsection{Two Lemmas}
\label{twolemmas}
In this section, we list some technical results. For $\bm{\Sigma}$
defined in \eqref{eq:improper_eta_sigma}, note that
 \begin{equation}
 \label{eq:improper_eta_sigma1}
  \bm{\Sigma}^{-1}=\left(\begin{array}{cc}
\left(\bm{X}^{T}\bm{X}\right)^{-1}+RS(\bm{\tau})^{-1}R^{T} & -RS(\bm{\tau})^{-1}\\
-S(\bm{\tau})^{-1}R^{T} & S(\bm{\tau})^{-1}
\end{array}\right),
 \end{equation}
with $S(\bm{\tau})$ and $R$ defined as
\begin{equation}
\label{eq:improper_eta_sigma2}
S(\bm{\tau})=\bm{Z}^{T}\left(\bm{I}-P_X\right)\bm{Z}+\bm{D}(\bm{\tau}), \text{ and }  R=\left(\bm{X}^{T}\bm{X}\right)^{-1}\bm{X}^{T}\bm{Z} 
\end{equation}
respectively, where 
\begin{equation}
\label{eq:projectionmatrix}
P_{X}=\bm{X}\left(\bm{X}^{T}\bm{X}\right)^{-1}\bm{X}^{T}.
\end{equation}
Also the mean for the conditional distribution of $\bm{\eta}$ in \eqref{eq:eta} becomes

\[\bm{\Sigma}^{-1}\bm{W}^{T}\bm{v}=\left(\begin{array}{c}
\left(\bm{X}^{T}\bm{X}\right)^{-1}\bm{X}^{T}\left[\bm{I}-\bm{Z}S(\bm{\tau})^{-1}\bm{Z}^{T}\left(\bm{I}-P_{X}\right)\right]\bm{v}\\
S(\bm{\tau})^{-1}\bm{Z}^{T}\left(\bm{I}-P_{X}\right)\bm{v}
\end{array}\right).\]

Let $U^{T}\Lambda U$ be the spectral decomposition of
$\bm{Z}^{T}\left(\bm{I}-P_{X}\right)\bm{Z}$ and let $\lambda_j$'s be
the diagonal elements of $\Lambda$.  Then
$\left(\bm{Z}^{T}\left(\bm{I}-P_{X}\right)\bm{Z}\right)^{+}\equiv
U^{T}\Lambda^{+}U$, where $\Lambda^{+}$ is a diagonal matrix whose
$j$th diagonal element is $\lambda_{j}^{+}=1/\lambda_{j}$ if
$\lambda_{j}\neq0$, and 0 otherwise.

\begin{lemma}
\label{probit_lemma1}
For the matrices $S(\bm{\tau})$ and $P_{X}$ defined in \eqref{eq:improper_eta_sigma2} and \eqref{eq:projectionmatrix}, the following inequalities hold for all $\tau_j \in \mathbb{R}_{+}$, $j=1,\dots,r$:
\begin{enumerate}
\item $S(\bm{\tau})^{-1}\preceq\left(\bm{Z}^{T}\left(\bm{I}-P_{X}\right)\bm{Z}\right)^{+}+\sum_{j=1}^{r}1/\tau_{j}\left(\bm{I}-P_{\bm{Z}^{T}\left(\bm{I}-P_{X}\right)\bm{Z}}\right)$.
\item $\left(R_{j}S(\bm{\tau})^{-1}R_{j}^{T}\right)^{-1}\preceq\left(\lambda_{p}+\tau_{j}\right)\bm{I}_{q_{j}}$,
where $\lambda_{p}$ is the largest eigenvalue of $\bm{Z}^{T}\left(\bm{I}-P_{X}\right)\bm{Z}$ and  $R_{j}$ is a $q_{j}\times q$ matrix with 0's and 1's such that $R_{j}\bm{u}=\bm{u}_{j}$.
\end{enumerate}

\end{lemma}
The proof of the above result is similar to that of Lemma 1 in \cite{roman2012convergence} and we omit it.
\begin{lemma}
\label{probit_lemma2}
Let $S(\bm{\tau})$ and $P_X$ be the two matrices as defined in \eqref{eq:improper_eta_sigma2} and \eqref{eq:projectionmatrix}. Let $\bm{l} = (l_1,\dots, l_n)^T \in \mathbb{R}^n$. For any $\bm{\tau} \in \mathbb{R}_+^r$, we have,
\[
\left\Vert S(\bm{\tau})^{-1}\bm{Z}^{T}\left(\bm{I-P}_{X}\right)\bm{l}\right\Vert \leq\hat{\varphi}\sum_{i=1}^{n}\left|l_{i}\right|,
\]
where $\Vert \cdot \Vert$ denotes the Euclidean norm and $\hat{\varphi}$ is a finite number that depends on $\bm{W}$.
\begin{proof}
Let $\bm{Z}_{P} \equiv\left(\bm{I}-P_{X}\right)\bm{Z}$ and $\bm{z}_{Pi}^{T}$
be the $i$th row of $\bm{Z}_{P}$. Then 
\begin{align*}
&\left\Vert S(\bm{\tau})^{-1}\bm{Z}^{T}\left(\bm{I-P}_{X}\right)\bm{l}\right\Vert   =\left\Vert \left(\bm{Z}_{P}^{T}\bm{Z}_{P}+\bm{D}(\bm{\tau})\right)^{-1}\bm{Z}_{P}^{T}\bm{l}\right\Vert \\
 &=\left\Vert \sum_{i=1}^{n}\left(\bm{Z}_{P}^{T}\bm{Z}_{P}+\bm{D}(\bm{\tau})\right)^{-1}\bm{z}_{Pi}l_{i}\right\Vert \leq\sum_{i=1}^{n}\left\Vert \left(\bm{Z}_{P}^{T}\bm{Z}_{P}+\bm{D}(\bm{\tau})\right)^{-1}\bm{z}_{Pi}l_{i}\right\Vert \\
 & =\sum_{i=1}^{n}\left\Vert \left(\sum_{k=1}^{n}\bm{z}_{Pk}\bm{z}_{Pk}^{T}+\bm{D}(\bm{\tau})\right)^{-1}\bm{z}_{Pi}l_{i}\right\Vert  \leq\sum_{i=1}^{n}\left|l_{i}\right| \varphi_{i}\left(\bm{\tau}\right),
\end{align*}
where 
\begin{equation*}
\varphi_{i}^{2}\left(\bm{\tau}\right) =\bm{z}_{Pi}^{T}\left(\bm{z}_{Pi}\bm{z}_{Pi}^{T}+\sum_{k\in\left\{ 1,\dots,n\right\} \backslash\left\{ i\right\} }\bm{z}_{Pk}\bm{z}_{Pk}^{T}+\bm{D}(\bm{\tau})\right)^{-2}\bm{z}_{Pi}.
\end{equation*}
Note that for fixed $i\in \{1,2\dots, n\}$,
\begin{align*}
 \varphi_{i}^{2}\left(\bm{\tau}\right) & =\bm{z}_{Pi}^{T}\left(\bm{z}_{Pi}\bm{z}_{Pi}^{T}+\sum_{k\in\left\{ 1,\dots,n\right\} \backslash\left\{ i\right\} }\bm{z}_{Pk}\bm{z}_{Pk}^{T}+\bm{D}(\bm{\tau}) \right. \\
 & \left.-\frac{1}{\sum_{j=1}^{r}1/\tau_{j}}\bm{I}_{q}+\frac{1}{\sum_{j=1}^{r}1/\tau_{j}}\bm{I}_{q}\right)^{-2}\bm{z}_{Pi}\\
 & \leq\sup_{\bm{\iota}\in \mathbb{R}_{+}^{n+q}}t_{i}^{T}\left(t_{i}t_{i}^{T}+\sum_{k\in\left\{ 1,\dots n\right\} \backslash\left\{ i\right\} }\iota_{k}t_{k}t_{k}^{T}+\sum_{k=n+1}^{n+q}\iota_{k}t_{k}t_{k}^{T}+\iota_{1}\bm{I}_{q}\right)^{-2}t_{i}\\
 & \equiv\hat{\varphi}_{i}^{2}, 
\end{align*}
where $\bm{\iota} = (\iota_1,\iota_2,\dots, \iota_{n+q})$,
$t_k = \bm{z}_{Pk}$ for $k \in \{1,\dots, n\}$ and for
$k = n+1,\dots, n+q$, define $t_k$ to be a $q\times 1$ unit vector
with 1 on the $(k-n)$th position, 0 elsewhere. The inequality
follows from the fact that $\sum_{j=1}^r 1/\tau_j > 1/\tau_j$.  By
Lemma 3 in \cite{roman2012convergence}, we know that
$\hat{\varphi}_{i}^{2}$ is finite.  Let
$\hat{\varphi}=\max_{1\leq i \leq n}\hat{\varphi}_{i}$, then
\[
\left\Vert S(\bm{\tau})^{-1}\bm{Z}^{T}\left(\bm{I-P}_{X}\right)\bm{l}\right\Vert \leq\sum_{i=1}^{n}\left|l_{i}\right|\hat{\varphi}_{i}\leq\hat{\varphi}\sum_{i=1}^{n}\left|l_{i}\right|.
\]
\end{proof}
\end{lemma}

\subsection{Proof of Theorem \ref{them2}}
\label{proof_thm2}

The two-block Gibbs sampler
$\{\bm{\eta}^{(m)}, (\bm{v}^{(m)}, \bm{\tau}^{(m)}) \}_{m=0}^{\infty}$
in Algorithm \ref{algorithm2} has the same rate of convergence as its
two marginal chains, namely, the $\bm{\eta}$-chain and the
$(\bm{v},\bm{\tau})$-chain. Here we work with the $\bm{\eta}$-chain, denoted as
$\bm{\Psi} = \{ \bm{\eta}^{(m)}\}_{m=0}^{\infty}$ and establish its
geometric rate of convergence. Define
$A \equiv \{ j\in \{1,\dots,r\}: b_j = 0\}$. Recall that given
$\bm{\eta}$, the conditional distribution of $\bm{\tau}$ is given by independent
$\text{Gamma} (a_j + q_j/2 , b_j + \bm{u}_j^T\bm{u}_j/2)$, $j=1,\dots,r$, which is
not defined when $A$ is not empty and
$\bm{\eta} \in \mathcal{N}=\left\{ \bm{\eta}\in\mathbb{R}^{p+q};
  \prod_{j\in A} ||\bm{u}_j|| =0\right\}$. Since $\mathcal{N}$ is a
set of measure zero, simulation of the Gibbs sampler is not affected
by the fact that $\pi(\bm{\tau}\vert \bm{\eta}, \bm{y})$ is not
defined on $\mathcal{N}$. But as mentioned in
\cite{roman2012convergence}, for a theoretical analysis of the
$\bm{\eta}$-chain, the Mtd of $\bm{\Psi}$ and hence
$\pi(\bm{\tau}\vert \bm{\eta}, \bm{y})$ must be defined for {\it all}
$\bm{\eta} \in \mathbb{R}^{p+q}$. Since $\mathcal{N}$ is a measure
zero set, the Mtd of $\bm{\Psi}$ hence
$\pi(\bm{\tau}\vert \bm{\eta}, \bm{y})$ can be defined arbitrarily on
$\mathcal{N}$. If $A$ is not empty for all
$\bm{\eta} \in \mathbb{R}^{p+q}$, we define
$\pi(\bm{\tau}\vert \bm{\eta}, \bm{y})$ as follows,
\[
\pi\left(\bm{\tau}|\bm{\eta},\bm{y}\right)=\begin{cases}
\prod_{j=1}^{r}f_{G}\left(\tau_{j},\frac{q_{j}}{2}+a_{j},\frac{\bm{u}_{j}^{T}\bm{u}_{j}}{2}+b_{j}\right) & \text{if }\bm{\eta}\notin\mathcal{N}\\
\prod_{j=1}^r f_{G}\left(\tau_{j},1,1\right) & \text{if }\bm{\eta}\in\mathcal{N}
\end{cases},
\]
where $f_G$ stands for the density of a Gamma random variable.

We denote the $\{\bm{\eta}^{(m)}\}_{m=0}^{\infty}$ Markov chain defined on $\mathbb{R}^{p+q} \backslash \mathcal{N}$ as $\tilde{\bm{\Psi}}$. The chain $\tilde{\bm{\Psi}}$ is Harris ergodic on $\mathbb{R}^{p+q} \backslash \mathcal{N}$. Our proof of geometric ergodicity of $\bm{\Psi}$  is through that of $\tilde{\bm{\Psi}}$. The following proof establishes the geometric ergodicity of $\tilde{\bm{\Psi}}$.

\begin{proof}
We prove the geometric ergodicity of $\tilde{\bm{\Psi}}$ by establishing a drift function, which has the following form,
\begin{equation}
\label{eq:drift1}
V\left(\bm{\eta}\right)=\sum_{i=1}^{n}\left(\bm{x}_{i}^{T}\bm{\beta}+\bm{z}_{i}^{T}\bm{u}\right)^{2}+\sum_{j=1}^{r}\left(\bm{u}_{j}^{T}\bm{u}_{j}\right)^{-c},
\end{equation}
where $c \in (0,1/2)$ is a positive constant determined later in the proof. Note that, since the condition A1 is in force, $V\left(\bm{\eta}\right) : \mathbb{R}^{p+q}\backslash\mathcal{N} \rightarrow [0,\infty)$ is unbounded off compact sets. 
We show that for any $\bm{\eta}, \bm{\eta}^{\prime} \in \mathbb{R}^{p+q} \backslash \mathcal{N}$, there exists $\rho_1 \in [0,1)$ and $L_1 >0$ such that
\begin{equation}
\label{eq:drift_ineq}
E[V\left(\bm{\eta}\right)|\bm{\eta}^{\prime}] \leq \rho_1V(\bm{\eta}^{\prime}) + L_1.
\end{equation}

By Fubini's theorem, we have
\begin{eqnarray*}
  E\left[V\left(\bm{\eta}\right)|\bm{\eta}^{\prime}\right] & = & \int_{\mathbb{R}^{p+q}\backslash \mathcal{N}} V\left(\bm{\eta}\right)k\left(\bm{\eta}|\bm{\eta}^{\prime}\right)d\bm{\eta}\\
                                                           & = & \int_{\mathbb{R}^n} \int_{\mathbb{R}^{r}_+} \int_{\mathbb{R}^{p+q}\backslash \mathcal{N}} V\left(\bm{\eta}\right)\pi\left(\bm{\eta}|\bm{v},\bm{\tau},\bm{y}\right)\pi\left(\bm{v},\bm{\tau}|\bm{\eta}^{\prime},\bm{y}\right)d\bm{\eta} d\bm{\tau} d\bm{v}.
\end{eqnarray*}
Thus, the expectation on the left hand side of \eqref{eq:drift_ineq} can be evaluated using two steps. First, we calculate the expectation with respect to the conditional distribution of $\bm{\eta}$ given $\bm{v}$, $\bm{\tau}$ and $\bm{y}$, that is $E[V(\bm{\eta})|\bm{v}, \bm{\tau},\bm{y}]$. 

From \eqref{eq:improper_eta_sigma1} and \eqref{eq:improper_eta_sigma2}, we have $\bm{W}\bm{\Sigma}^{-1}\bm{W}^{T}=P_{X}+\left(\bm{I}-P_{X}\right)\bm{Z}S(\bm{\tau})^{-1}\bm{Z}^{T}\left(\bm{I}-P_{X}\right)$.  Also $\left(\bm{I}-P_{X}\right)=\left(\bm{I}-P_{X}\right)^{2}$. Let $\tilde{P}=\left(\bm{I}-P_{X}\right)\bm{Z}\bm{D}(\bm{\tau})^{-1/2}$, then
\begin{align*}
& \left(\bm{I}-P_{X}\right)\bm{Z}S(\bm{\tau})^{-1}\bm{Z}^{T}\left(\bm{I}-P_{X}\right) \\
=& \left(\bm{I}-P_{X}\right)^{2}\bm{Z}S(\bm{\tau})^{-1}\bm{Z}^{T}\left(\bm{I}-P_{X}\right)^{2}\\
= & \left(\bm{I}-P_{X}\right)\tilde{P}\left(\tilde{P}^{T}\tilde{P}+\bm{I}\right)^{-1}\tilde{P}^{T}\left(\bm{I}-P_{X}\right)\preceq \bm{I}-P_{X}.
\end{align*}
 Thus, $\bm{W}\bm{\Sigma}^{-1}\bm{W}^{T}\preceq P_{X}+\bm{I}-P_{X}=\bm{I}$.  Here ``$\bm{W}\bm{\Sigma}^{-1}\bm{W}^{T} \preceq \bm{I}$"  means that $\bm{I} - \bm{W}\bm{\Sigma}^{-1}\bm{W}^{T} $ is a positive semidefinite matrix.  From \eqref{eq:eta} and \eqref{eq:improper_eta_sigma}, it follows that 
\begin{align}
&E\left[\sum_{i=1}^{n}\left(\bm{x}_{i}^{T}\bm{\beta}+\bm{z}_{i}^{T}\bm{u}\right)^{2}|\bm{v},\bm{\tau}, \bm{y}\right]  \leq E\left[\bm{\eta}^{T}\bm{\Sigma}\bm{\eta}|\bm{\tau},\bm{v}, \bm{y}\right] \nonumber\\
 & =p+q+\bm{v}\bm{W}\bm{\Sigma}^{-1}\bm{W}^{{T}}\bm{v} \leq p+q+\bm{v}^T\bm{v}. \label{eq:improper_thm2_exp1}
\end{align}

According to \cite{roman2012convergence}, for $c\in (0,1/2)$ we have
\begin{equation}
E\left[\left(\bm{u}_{j}^{T}\bm{u}_{j}\right)^{-c}|\bm{v},\bm{\tau}, \bm{y}\right]\leq2^{-c}\frac{\Gamma\left(q_{j}/2-c\right)}{\Gamma\left(q_{j}/2\right)}\left[\lambda_{p}^{c}+\tau_{j}^{ c}\right], 
\label{eq:improper_thm2_exp2}
\end{equation}
where $\lambda_p$ is the largest eigenvalue of $\bm{Z}^T(\bm{I}-P_{\bm{X}})\bm{Z}$. Using \eqref{eq:improper_thm2_exp1} and \eqref{eq:improper_thm2_exp2} from \eqref{eq:drift1}, we have
\begin{equation}
\label{eq:improper_thm2_step1}
E\left[V\left(\bm{\eta}\right)|\bm{\tau},\bm{v},\bm{y}\right]\leq\ \bm{v}^{T}\bm{v}+2^{-c}\sum_{j=1}^{r}\frac{\Gamma\left(q_{j}/2-c\right)}{\Gamma\left(q_{j}/2\right)}\tau_{j}^{c} + 2^{-c}\sum_{j=1}^{r}\frac{\Gamma\left(q_{j}/2-c\right)}{\Gamma\left(q_{j}/2\right)}\lambda_p^c+p+q.
\end{equation}

Now we consider the expectation corresponding to the conditional distribution of $\bm{v}$ and $\bm{\tau}$ given $\bm{\eta}^{\prime}$ and $\bm{y}$. Using (10) from \cite{roy2007convergence}, we have
\[
E(v_{i}^{2}|\bm{\eta}^{\prime},\bm{y})=\left\{ \begin{array}{ccc}
1+(\bm{x}_{i}^{T}\bm{\beta}^{\prime}+\bm{z}_{i}^{T}\bm{u}^{\prime})^{2}+\frac{\left(\bm{x}_{i}^{T}\bm{\beta}^{\prime}+\bm{z}_{i}^{T}\bm{u}^{\prime}\right)\phi\left(\bm{x}_{i}^{T}\bm{\beta}^{\prime}+\bm{z}_{i}^{T}\bm{u}^{\prime}\right)}{\Phi\left(\bm{x}_{i}^{T}\bm{\beta}^{\prime}+\bm{z}_{i}^{T}\bm{u}^{\prime}\right)} &  & \text{if }y_{i}=1\\
1+(\bm{x}_{i}^{T}\bm{\beta}^{\prime}+\bm{z}_{i}^{T}\bm{u}^{\prime})^{2}-\frac{\left(\bm{x}_{i}^{T}\bm{\beta}^{\prime}+\bm{z}_{i}^{T}\bm{u}^{\prime}\right)\phi\left(\bm{x}_{i}^{T}\bm{\beta}^{\prime}+\bm{z}_{i}^{T}\bm{u}^{\prime}\right)}{1- \Phi\left(\bm{x}_{i}^{T}\bm{\beta}^{\prime}+\bm{z}_{i}^{T}\bm{u}^{\prime}\right)} &  & \text{if }y_{i}=0
\end{array}\right. .
\]
The above expectation can be written as,
\begin{equation}
\label{eq:exp_v}
E(v_{i}^{2}|\bm{\eta}^{\prime},\bm{y})=1+\left(\bm{w}_{i}^{*T}\bm{\eta}^{\prime}\right)^{2}-\frac{\left(\bm{w}_{i}^{*T}\bm{\eta}^{\prime}\right)\phi\left(\bm{w}_{i}^{*T}\bm{\eta}^{\prime}\right)}{1-\Phi\left(\bm{w}_{i}^{*T}\bm{\eta}^{\prime}\right)},
\end{equation}
where $\bm{w}_i^* = c_i\bm{w}_i^T$ is the $i$th row of $\bm{W}^*$ defined in section \ref{sec_propriety}. Also,
\begin{eqnarray}
-\frac{\left(\bm{w}_{i}^{*T}\bm{\eta}^{\prime}\right)\phi\left(\bm{w}_{i}^{*T}\bm{\eta}^{\prime}\right)}{1-\Phi\left(\bm{w}_{i}^{*T}\bm{\eta}^{\prime}\right)} & \leq & \begin{cases}
\left|\frac{\left(\bm{w}_{i}^{*T}\bm{\eta}^{\prime}\right)\phi\left(\bm{w}_{i}^{*T}\bm{\eta}^{\prime}\right)}{1-\Phi\left(\bm{w}_{i}^{*T}\bm{\eta}^{\prime}\right)}\right| & \text{if }\bm{w}_{i}^{*T}\bm{\eta}^{\prime}\leq0\\
0 & \text{if }\bm{w}_{i}^{*T}\bm{\eta}^{\prime}>0
\end{cases}\nonumber\\
 & \leq & \sup_{u\in\left(-\infty,0\right]}\left|\frac{u\phi\left(u\right)}{1-\Phi\left(u\right)}\right|\equiv\Xi, \label{eq:vbound}
\end{eqnarray}
where $\Xi \in (0,\infty)$. 

We use $A_{1},\dots,A_{2^{n}}$ to denote all the subsets of $\mathbb{N}_{n}=\left\{ 1,2\dots,n\right\} $. Following \cite{roy2007convergence}, let 
\[
S_{j}=\left\{ \bm{\eta}^{\prime}\in\mathbb{R}^{p+q}\backslash\left\{ \bm{0}\right\} :\bm{w}_{i}^{T}\bm{\eta}^{\prime}\leq0\text{ for all }i\in A_{j}\text{ and }\bm{w}_{i}^{T}\bm{\eta}^{\prime}>0\text{ for all }i\in\bar{A}_{j}\right\},
\]
where $\bar{A}_{j}$ is the complement of $A_{j}$. As mentioned in \cite{roy2007convergence}, the sets $S_j$'s are disjoint, $\cup_{j=1}^{2^n}S_j = \mathbb{R}^{p+q}\backslash \{\bm{0}\}$ and some of the $S_j$'s may be empty. For $j\in C\equiv\left\{ i\in\mathbb{N}_{2^{n}}:S_{i}\neq\emptyset\right\} $, define 
\[
H_{j}\left(\bm{\eta}^{\prime}\right)=\frac{\sum_{i\in A_{j}}\left(\bm{w}_{i}^{*T}\bm{\eta}^{\prime}\right)^{2}}{\sum_{i=1}^{n}\left(\bm{w}_{i}^{*T}\bm{\eta}^{\prime}\right)^{2}}=\frac{\sum_{i\in A_{j}}\left(\bm{w}_{i}^{*T}\bm{\eta}^{\prime}\right)^{2}}{\sum_{i\in A_{j}}\left(\bm{w}_{i}^{*T}\bm{\eta}^{\prime}\right)^{2}+\sum_{i\in\bar{A_{j}}}\left(\bm{w}_{i}^{*T}\bm{\eta}^{\prime}\right)^{2}}.
\]
 By \eqref{eq:exp_v}, for $\bm{\eta}\in S_{j}$, $j\in C$, we have 
\begin{align*}
 &E\left[\sum_{i=1}^{n}v_{i}^{2}|\bm{\eta}^{\prime}, \bm{y}\right] \\
 =& n+\sum_{i=1}^{n}\left(\bm{w}_{i}^{*T}\bm{\eta}^{\prime}\right)^{2}-\sum_{i\in A_{j}}\frac{\left(\bm{w}_{i}^{*T}\bm{\eta}^{\prime}\right)\phi\left(\bm{w}_{i}^{*T}\bm{\eta}^{\prime}\right)}{1-\Phi\left(\bm{w}_{i}^{*T}\bm{\eta}^{\prime}\right)}-\sum_{i\in\bar{A_{j}}}\frac{\left(\bm{w}_{i}^{*T}\bm{\eta}^{\prime}\right)\phi\left(\bm{w}_{i}^{*T}\bm{\eta}^{\prime}\right)}{1-\Phi\left(\bm{w}_{i}^{*T}\bm{\eta}^{\prime}\right)}\\
 =& n+\sum_{i=1}^{n}\left(\bm{w}_{i}^{*T}\bm{\eta}^{\prime}\right)^{2}+\sum_{i\in A_{j}}\left|\frac{\left(\bm{w}_{i}^{*T}\bm{\eta}^{\prime}\right)\phi\left(\bm{w}_{i}^{*T}\bm{\eta}^{\prime}\right)}{1-\Phi\left(\bm{w}_{i}^{*T}\bm{\eta}^{\prime}\right)}\right|-\sum_{i\in\bar{A_{j}}}\frac{\left(\bm{w}_{i}^{*T}\bm{\eta}^{\prime}\right)\phi\left(\bm{w}_{i}^{*T}\bm{\eta}^{\prime}\right)}{1-\Phi\left(\bm{w}_{i}^{*T}\bm{\eta}^{\prime}\right)}\\
 \leq &  n+\sum_{i=1}^{n}\left(\bm{w}_{i}^{*T}\bm{\eta}^{\prime}\right)^{2}+n\Xi-\sum_{i\in\bar{A_{j}}}\left(\bm{w}_{i}^{*T}\bm{\eta}^{\prime}\right)^{2}\\
 = & n\left(1+\Xi\right)+H_{j}\left(\bm{\eta}^{\prime}\right)\sum_{i=1}^{n}\left(\bm{w}_{i}^{*T}\bm{\eta}^{\prime}\right)^{2},
\end{align*}
where $\Xi$ is defined in \eqref{eq:vbound} and the inequality is due to the fact that $u\phi\left(u\right)/\left[1-\Phi\left(u\right)\right]\geq u^{2}$ for $u\geq0$.  Define $\lambda_{j}=\sup_{\bm{\eta}^{\prime}\in S_{j}}\left\{ H_{j}\left(\bm{\eta}^{\prime}\right)\right\} \in\left[0,1\right]$
and 
\[
\lambda_{0}=\max_{j\in C}\lambda_{j}.
\] 

If $\bm{\eta}^{\prime}=\bm{0}$, from \eqref{eq:exp_v}, we have $E\left[\sum_{i=1}^{n}v_{i}^{2}|\bm{\eta}^{\prime}, \bm{y}\right]=n$. Thus, for all $\bm{\eta}^{\prime} \in \mathbb{R}^{p+q}$, 
\begin{equation}
\label{eq:improper_thm2_ev}
E\left[\sum_{i=1}^{n}v_{i}^{2}|\bm{\eta}^{\prime}, \bm{y}\right]\leq\lambda_{0}\sum_{i=1}^{n}\left(\bm{x}_{i}^{T}\bm{\beta}^{\prime}+\bm{z}_{i}^{T}\bm{u}^{\prime}\right)^{2}+n\left(1+\Xi\right).
\end{equation}
Since conditions A1 and A2 are in force, using the techniques in \cite{roy2007convergence}, it can be shown that $\lambda_{0}<1$.

For $c\in (0, 1/2)$, define 
\begin{equation}
\label{eq:gamma_g}
G_j(-c) = 2^c \frac{\Gamma(q_j/2+a_j+c)}{\Gamma(q_j/2 + a_j)} \text{ for } j=1,\dots, r.
\end{equation}
Since $\tau_j\vert \bm{\eta}^{\prime}, \bm{y} \sim \text{Gamma} (a_j + q_j/2, b_j + \bm{u_j}^{\prime T}\bm{u}_j^{\prime}/2)$, 
\begin{align}
\label{eq:improper_thm2_tau}
E\left[\tau_{j}^{c}|\bm{u}_{j}^{\prime},\bm{y}\right]&=2^{-c}G_{j}\left(-c\right)\left[b_{j}+\frac{\bm{u}_{j}^{\prime T}\bm{u}_{j}^{\prime }}{2}\right]^{-c} \nonumber\\
 &\leq   G_{j}\left(-c\right)\left[\left(2b_{j}\right)^{-c}I_{\left(0,\infty\right)}\left(b_{j}\right)+\left(\bm{u}_{j}^{\prime T}\bm{u}_{j}^{\prime}\right)^{-c}I_{\left\{ 0\right\} }\left(b_{j}\right)\right].
\end{align}

Recall that $A = \{j\in \{1, 2,\dots, r\}: b_j = 0 \}$. We consider two cases, namely, when $A$ is empty and $A$ is not empty. 

\textbf{Case 1}: $A$ is not empty.

Then using \eqref{eq:improper_thm2_ev} and \eqref{eq:improper_thm2_tau}, from \eqref{eq:improper_thm2_step1} we have
\[
E\left[V\left(\bm{\eta}\right)|\bm{\eta}^{\prime}\right]\leq \lambda_0\sum_{i=1}^{n}\left(\bm{x}_{i}^{T}\bm{\beta}^{\prime}+\bm{z}_{i}^{T}\bm{u}^{\prime}\right)^2+\delta_{1}\left(c\right)\sum_{j\in A}\left(\bm{u}_{j}^{\prime T}\bm{u}_{j}^{\prime}\right)^{-c}+L_1\left(c\right),
\]
where 
\begin{align}
\delta_{1}\left(c\right) & \equiv 2^{-c}\max_{j\in A}G_{j}\left(-c\right)\frac{\Gamma\left(q_{j}/2-c\right)}{\Gamma\left(q_{j}/2\right)}, \label{eq:delta1c}\\
L_1\left(c\right) &  \equiv n(1+\Xi)+p+q+2^{-c}\lambda_{p}^{c}\sum_{j=1}^{r}\frac{\Gamma\left(q_{j}/2-c\right)}{\Gamma\left(q_{j}/2\right)} \nonumber \\
&+2^{-c}\sum_{j\notin A}G_{j}\left(-c\right)\frac{\Gamma\left(q_{j}/2-c\right)}{\Gamma\left(q_{j}/2\right)}\left(2b_{j}\right)^{-c}. \nonumber 
\end{align}

By \cite{roman2012convergence}, there exists $c\in C_1 = (0,1/2)\cap(0,-\max_{j\in A}a_j) $ such that $\delta_1(c)<1$. Thus, taking $\rho_1 = \max(\lambda_0,\delta_1(c))$, and $L_1 = L_1(c)$, we have

\[
E\left[V\left(\bm{\eta}\right)|\bm{\eta}^{\prime}\right]\leq  \rho_1 V\left( \bm{\eta}^{\prime} \right)  + L_1.
\]

\textbf{Case 2}: $A$ is empty.

In this case, the conditional expectation of $\tau_j^c$ can be bounded by a constant. Indeed from \eqref{eq:improper_thm2_tau} we have
\[
E\left[\tau_{j}^{c}|\bm{u}_{j}^{\prime}, \bm{y}\right]=2^{-c}G_{j}\left(-c\right)\left[b_{j}+\frac{\bm{u}_{j}^{\prime T}\bm{u}_{j}^{\prime }}{2}\right]^{-c}\leq G_{j}\left(-c\right)\left(2b_{j}\right)^{-c}.
\]
Thus when $A$ is empty, we have
\[
E\left[V\left(\bm{\eta} \right)|\bm{\eta}^{\prime}\right]\leq \lambda_0\sum_{i=1}^{n}\left(\bm{x}_{i}^{T}\bm{\beta}^{\prime}+\bm{z}_{i}^{T}\bm{u}^{\prime}\right)^2+L_1(c)
\leq  \lambda_0V\left( \bm{\eta}^{\prime} \right)  + L_1(c).
\]

Hence in both cases, \eqref{eq:drift_ineq} holds.  We now show that $\bm{\eta}$-chain is a Feller chain on $\mathbb{R}^{p+q} \backslash \mathcal{N}$, which means that $K\left(\bm{\eta},O\right)$ is a lower semi-continuous function on $\mathbb{R}^{p+q} \backslash \mathcal{N}$ for each fixed open set $O$. For a sequence $\{\bm{\eta}_m\}$ note that,
\begin{align*}
\liminf_{m\rightarrow\infty} K\left(\bm{\eta}_{m},O\right) & =\liminf_{m\rightarrow\infty}\int_{O}k\left(\bm{\eta}|\bm{\eta}_{m}\right)d\bm{\eta}\\
 & =\liminf_{m\rightarrow\infty}\int_{O}\left[\int_{\mathbb{R}_{+}^{r}}\int_{\mathbb{R}^{n}}\pi(\bm{\eta}|\bm{v},\bm{\tau},\bm{y})\pi(\bm{v},\bm{\tau}|\bm{\eta}_{m},\bm{y})d\bm{v}d\bm{\tau}\right]d\bm{\eta}\\
 & \geq\int_{O}\int_{\mathbb{R}_{+}^{r}}\int_{\mathbb{R}^{n}}\pi(\bm{\eta}|\bm{v},\bm{\tau},\bm{y})\liminf_{m\rightarrow\infty}\pi(\bm{v},\bm{\tau}|\bm{\eta}_{m},\bm{y})d\bm{v}d\bm{\tau}d\bm{\eta},
\end{align*}
where the inequality follows from Fatou's lemma. Recall that $\pi(\bm{v},\bm{\tau}|\bm{\eta},\bm{y})=\pi(\bm{v}|\bm{\eta},\bm{y})\pi(\bm{\tau}|\bm{\eta},\bm{y})$. Note that, $\tau_j \vert \bm{\eta}^{\prime}, \bm{y} \sim \text{Gamma}(a_j + q_j/2, b_j + \bm{u}_j^{\prime T}\bm{u}_j^{\prime}/2)$ and condition A3 holds. Thus, for all $\bm{\eta}^{\prime}\in\mathbb{R}^{p+q} \backslash \mathcal{N}$ the conditional distribution of $\tau_j$ is a Gamma distribution with positive shape and scale parameters even if $b_j = 0$. Since both $\pi(\bm{v}|\bm{\eta},\bm{y})$ and $\pi(\bm{\tau}|\bm{\eta},\bm{y})$ are continuous functions in $\bm{\eta} \in \mathbb{R}^{p+q} \backslash \mathcal{N}$ ,  if $\bm{\eta}_m\rightarrow \bm{\eta}$,
\begin{align*}
\liminf_{m\rightarrow\infty} K\left(\bm{\eta}_{m},O\right) & \geq\int_{O}\int_{\mathbb{R}_{+}^{r}}\int_{\mathbb{R}^{n}}\pi(\bm{\eta}|\bm{v},\bm{\tau},\bm{y})\pi(\bm{v},\bm{\tau}|\bm{\eta},\bm{y})d\bm{v}d\bm{\tau}d\bm{\eta}\\
 & =K\left(\bm{\eta},O\right).
\end{align*}
Thus by \cite{meyn1993markov}(chap. 15), \eqref{eq:drift_ineq} implies the Markov chain $\tilde{\bm{\Psi}}$ is geometrically ergodic.

Next, we need to show that the original Markov chain $\bm{\Psi}$ is geometrically ergodic. The techniques of Lemma 12 in \cite{roman2012thesis} can be applied here for this purpose. 

Let $M$ and $\tilde{M}$ be the Mtfs of $\bm{\Psi}$ and
$\tilde{\bm{\Psi}}$ respectively. Also, let $M^m$ and
$\tilde{M}^m$ be the corresponding $m$-step Mtfs, and
$\mathsf{X} \equiv \mathbb{R}^{p+q}$,
$\tilde{\mathsf{X}} \equiv \mathbb{R}^{p+q} \backslash \mathcal{N}$.
Recall that $\mathscr{B}$ denotes the Borel $\sigma$-algebra of
$\mathbb{R}^{p+q}$.  Since the Lebesgue measure of $\mathcal{N}$ is 0,
for any $\mathsf{x} \in \tilde{\mathsf{X}}$ and
$\mathsf{B} \in \mathcal{B}_{\tilde{\mathsf{X}}} = \{
\tilde{\mathsf{X}} \cap \mathsf{A}: \mathsf{A} \in \mathscr{B}\}$
\[
\tilde{M}(\mathsf{x}, \mathsf{B}) = M(\mathsf{x}, \mathsf{B}).
\]

Let $\mu$ and $\tilde{\mu}$ be the Lebesgue measures on $\mathsf{X}$
and $\tilde{\mathsf{X}}$ respectively. Then $\bm{\Psi}$ and
$\tilde{\bm{\Psi}}$ are $\mu$-irreducible and
$\tilde{\mu}$-irreducible respectively. Also, $\mu$ and $\tilde{\mu}$
are the corresponding maximal irreducibility measures. These two
Markov chains $\bm{\Psi}$ and $\tilde{\bm{\Psi}}$ are also
aperiodic. According to Theorem 15.0.1 in \cite{meyn1993markov}, there
exists a $\nu$-petite set
$\mathsf{C} \in \mathcal{B}_{\tilde{\mathsf{X}}} $,
$\rho_{\mathsf{C}} < 1$, $M_{\mathsf{C}} < \infty$, a number
$\tilde{M}^{\infty}(\mathsf{C})$ such that
$\tilde{\mu}(\mathsf{C}) >0$ and
\begin{equation*}
|\tilde{M}^m(\mathsf{x},\mathsf{C}) - \tilde{M}^{\infty}(\mathsf{C})|<M_{\mathsf{C}}\rho^m_{\mathsf{C}},
\end{equation*}
for all $\mathsf{x} \in \mathsf{C}$. Since the set $\mathsf{C}$ is a $\nu$-petite set for $\tilde{\mathsf{X}}$,  $\nu$ is a nontrivial measure on $\mathcal{B}_{\tilde{\mathsf{X}}}$ with, 
\begin{equation*}
\sum_{m=0}^{\infty}\tilde{M}^m(\mathsf{x},\mathsf{B})\tilde{a}(m) \geq \nu(\mathsf{B})
\end{equation*}
for all $\mathsf{x} \in \mathsf{C}$ and $\mathsf{B} \in \mathcal{B}_{\tilde{\mathsf{X}}}$, where $\tilde{a}(m)$ is a mass function on $\{0,1,2,\dots,\}$. 

Since $\tilde{M}^m(\mathsf{x}, \mathsf{B}) = M^m(\mathsf{x}, \mathsf{B})$ for any $\mathsf{x} \in \tilde{\mathsf{X}}$ and $\mathsf{B} \in \mathcal{B}_{\tilde{\mathsf{X}}}$, we have $ M^m(\mathsf{x},\mathsf{C})=\tilde{M}^m(\mathsf{x},\mathsf{C})$. So for all $x \in \mathsf{C}$
\begin{equation*}
|M^m(\mathsf{x},\mathsf{C}) - \tilde{M}^{\infty}(\mathsf{C})|<M_{\mathsf{C}}\rho^m_{\mathsf{C}}.
\end{equation*}

Also, since $\mu(\mathcal{N}) = 0$,  we know that $\mu(\mathsf{C}) >0$. It can be checked that $\mathsf{C}$ is also petite for the original Markov chain $\bm{\Psi}$. Thus from Theorem 15.0.1 of \cite{meyn1993markov}, it follows that $\bm{\Psi}$ is geometrically ergodic.

\end{proof}

\subsection{Proof of Theorem \ref{them4}}
\label{proof_thm4}
\begin{proof}
As in Appendix \ref{proof_thm2}, we study the convergence properties of the $\bm{\eta}$-chain. Recall that $\mathcal{N}=\left\{ \bm{\eta}\in\mathbb{R}^{p+q}; \prod_{j\in A} ||\bm{u}_j|| =0\right\}$. When $A$ is nonempty and $\bm{\eta} \in \mathcal{N}$, we define the conditional distribution of $\bm{\tau}$ given $\bm{\eta}, \bm{y}$ the same way as in Appendix \ref{proof_thm2}.

Consider the following drift function on $\mathbb{R}^{p+q} \backslash \mathcal{N}$,
\[
V\left(\bm{\eta}\right)=\alpha\sum_{i=1}^{n}\left(\bm{x}_{i}^{T}\bm{\beta}+\bm{z}_{i}^{T}\bm{u}\right)^{2}+\sum_{j=1}^{r}G_{j}\left(s\right)\left(\bm{u}_{j}^{T}\bm{u}_{j}\right)^{s}+\sum_{j=1}^{r}\left(\bm{u}_{j}^{T}\bm{u}_{j}\right)^{-c}.
\]
where $G_j(\cdot)$ is defined in \eqref{eq:gamma_g}, $\alpha$,
$s \in \tilde{S} \equiv (0,1]\cap\left(0,\tilde{s}\right)$ for
$\tilde{s}$ defined in Theorem \ref{them4}, and
$c \in C_1 = (0,1/2)\cap(0,-\max_{j\in A}a_j) $ are positive constants
to be chosen later. Under the assumption B3,
$V(\bm{\eta}): \mathbb{R}^{p+q} \backslash \mathcal{N} \rightarrow
[0,\infty)$ is unbounded off compact sets (Since $W$ is not a full
rank matrix, the drift function considered in the proof of Theorem
\ref{them2} is no more unbounded off compact sets.). We need to show
that for any
$\bm{\eta}, \bm{\eta}^{\prime} \in \mathbb{R}^{p+q} \backslash
\mathcal{N}$, there exists a constant $\rho_2 \in [0,1)$ and $L_2>0$
such that
\begin{equation}
\label{drift4}
E[V\left(\bm{\eta}\right)|\bm{\eta}^{\prime}]  = E\{ E [V(\bm{\eta}\vert \bm{v, \tau, y})]\vert \bm{\eta}^{\prime}, \bm{y}\}\leq \rho_2V(\bm{\eta}^{\prime}) + L_2.
\end{equation}

First, we calculate the expectation of $V(\bm{\eta})$ with respect to the $\bm{\eta}$ conditional distribution given $\bm{v, \tau}$ and $\bm{y}$. Same calculations as in the proof of Theorem \ref{them2} (see \eqref{eq:improper_thm2_exp1}) show that,
\begin{equation}
\label{eq:improper_thm4_exp1}
E\left[\sum_{i=1}^{n}\left(\bm{x}_{i}^{T}\bm{\beta}+\bm{z}_{i}^{T}\bm{u}\right)^{2}\vert \bm{v}, \bm{\tau}, \bm{y}\right] 
 \leq p+q+\bm{v}^{T}\bm{v}.
\end{equation}

For $s\in(0,1]$, by Jensen inequality,
\begin{equation}
\label{eq:improper_thm4_exp2}
E\left[\left(\bm{u}_{j}^{T}\bm{u}_{j}\right)^{s}|\bm{v},\bm{\tau}, \bm{y}\right]\leq\left[E\left(\bm{u}_{j}^{T}\bm{u}_{j}|\bm{v},\bm{\tau}, \bm{y}\right)\right]^{s}.
\end{equation}
 Also, from \eqref{eq:eta} and \eqref{eq:improper_eta_sigma} it follows that 
\begin{equation}
\label{eq:improper_thm4_exp2ineq}
E\left(\bm{u}_{j}^{T}\bm{u}_{j}|\bm{v},\bm{\tau}, \bm{y}\right)=tr\left(R_{j}S(\bm{\tau})^{-1}R_{j}\right)+\left[E\left(R_{j}\bm{u}|\bm{v},\bm{\tau}, \bm{y}\right)\right]^{T}\left[E\left(R_{j}\bm{u}|\bm{v},\bm{\tau}, \bm{y}\right)\right],
\end{equation}
where $R_j$ is defined in Lemma \ref{probit_lemma1}. For the first part on the right hand side of \eqref{eq:improper_thm4_exp2ineq}, we have
\begin{align}
tr\left(R_{j}S(\bm{\tau})^{-1}R_{j}^T\right) & =tr\left[R_{j}\left(\bm{Z}^{T}\left(\bm{I}-P_{X}\right)\bm{Z}\right)^{+}R_{j}^{T}\right] \nonumber \\
&+tr\left[R_{j}\left(\bm{I}-P_{\bm{Z}^{T}\left(\bm{I}-P_{X}\right)\bm{Z}}\right)R_{j}^{T}\right]\sum_{l=1}^{r}\tau_{l}^{-1} \nonumber \\
 & =\xi_{j}+\varsigma_{j}\sum_{l=1}^{r}\tau_{l}^{-1},  \label{eq:improper_thm4_trace}
\end{align}
where $ \xi_{j} = tr\left[R_{j}\left(\bm{Z}^{T}\left(\bm{I}-P_{X}\right)\bm{Z}\right)^{+}R_{j}^{T}\right]$ and $\varsigma_{j} = tr\left[R_{j}\left(\bm{I}-P_{\bm{Z}^{T}\left(\bm{I}-P_{X}\right)\bm{Z}}\right)R_{j}^{T}\right]$. For the second part,  we have
\begin{align}
&\left[E\left(R_{j}\bm{u}|\bm{v},\bm{\tau}, \bm{y}\right)\right]^{T}\left[E\left(R_{j}\bm{u}|\bm{v},\bm{\tau}, \bm{y}\right)\right]  \nonumber \\
&=\bm{v}^{ T}\left(\bm{I}-P_{X}\right)\bm{Z}S(\bm{\tau})^{-1}R_{j}^{T}R_{j}S(\bm{\tau})^{-1}\bm{Z}^{T}\left(\bm{I}-P_{X}\right)\bm{v} \nonumber \\
 & \leq\bm{v}^{T}\left(\bm{I}-P_{X}\right)\bm{Z}S(\bm{\tau})^{-1}S(\bm{\tau})^{-1}\bm{Z}^{T}\left(\bm{I}-P_{X}\right)\bm{v} \nonumber \\
 & =\left\Vert S(\bm{\tau})^{-1}\bm{Z}^{T}\left(\bm{I}-P_{X}\right)\bm{v}\right\Vert ^{2} \nonumber \\
 & \leq\left(\hat{\varphi}\sum_{i=1}^{n}\left|v_{i}\right|\right)^{2}\leq\hat{\varphi}^2n\sum_{i=1}^{n}v_{i}^{2}, \label{eq:improper_thm4_s2}
\end{align}
where the second inequality follows from Lemma \ref{probit_lemma2} given in Appendix \ref{twolemmas}. Combining \eqref{eq:improper_thm4_trace} and \eqref{eq:improper_thm4_s2}, from \eqref{eq:improper_thm4_exp2ineq} we have
\begin{align*}
\left[E\left(\bm{u}_{j}^{T}\bm{u}_{j}|\bm{v},\bm{\tau}, \bm{y}\right)\right]^{s} &\leq\left[\xi_{j}+\varsigma_{j}\sum_{j=1}^{r}\tau_j^{ -1}+\hat{\varphi}^{2}n\sum_{i=1}^{n}v_{i}^{2}\right]^{s}\\
&\leq\xi_{j}^{s}+\varsigma_{j}^{s}\sum_{l=1}^{r}\tau_{l}^{-s}+\hat{\varphi}^{2s}n^{s}\sum_{i=1}^{n}v_{i}^{2s}.
\end{align*}

Note that, if $v_{i}^{2}\leq1$, then $v_{i}^{2s}\leq1$, 
and if $v_{i}^{2s}>1$, then $v_{i}^{2s}<v_{i}^{2}$. So $v_{i}^{2s}\leq1+v_{i}^{2}$. Thus, 
\begin{equation}
\label{eq:improper_thm4_ubound}
\left[E\left(\bm{u}_{j}^{T}\bm{u}_{j}|\bm{v},\bm{\tau}, \bm{y}\right)\right]^{s}\leq\varsigma_{j}^{s}\sum_{l=1}^{r}\tau_{l}^{-s}+\hat{\varphi}^{2s}n^{s}\sum_{i=1}^{n}v_{i}^{2}+\hat{\varphi}^{2s}n^{1+s}+\xi_{j}^{s}.
\end{equation}

Also recall from \eqref{eq:improper_thm2_exp2} that
we also have,
\[
E\left[\left(\bm{u}_{j}^{T}\bm{u}_{j}\right)^{-c}|\bm{v},\bm{\tau}, \bm{y}\right]\leq2^{-c}\frac{\Gamma\left(q_{j}/2-c\right)}{\Gamma\left(q_{j}/2\right)}\left[\lambda_{p}^{c}+\tau_{j}^{c}\right].
\] Combining \eqref{eq:improper_thm2_exp2} , \eqref{eq:improper_thm4_exp1}, \eqref{eq:improper_thm4_exp2} and \eqref{eq:improper_thm4_ubound} from \eqref{drift4} we have
\begin{align}
\label{eq:improper_thm4_step1}
E\left[V\left(\bm{\eta}\right)|\bm{v},\bm{\tau}, \bm{y}\right]&\leq\left(\alpha+\delta_{2}\left(s\right)\right)\sum_{i=1}^{n}v_{i}^{2}+\delta_{3}\left(s\right)\sum_{j=1}^{r}\tau_{j}^{-s} \nonumber\\
&+2^{-c}\sum_{j=1}^{r}\frac{\Gamma\left(q_{j}/2-c\right)}{\Gamma\left(q_{j}/2\right)}\tau_{j}^{c}+\kappa_{1}\left(\alpha,s,c\right),
\end{align}
where 
\begin{align*}
\delta_{2}\left(s\right) & =\hat{\varphi}^{2s}n^{s}\sum_{j=1}^{r}G_{j}\left(s\right),\\
\delta_{3}\left(s\right) & =\sum_{j=1}^{r}G_{j}\left(s\right)\varsigma_{j}^{s}, \,\text{and}\\
\kappa_{1}\left(\alpha,s,c\right) & =\alpha\left(p+q\right)+\sum_{j=1}^{r}G_{j}\left(s\right)\left(\hat{\varphi}^{2s}n^{1+s}+\xi_{j}^{s}\right)+2^{-c}\lambda_{p}^{c}\sum_{j=1}^{r}\frac{\Gamma\left(q_{j}/2-c\right)}{\Gamma\left(q_{j}/2\right)}.
\end{align*}

Next we calculate the outer expectation in \eqref{drift4}, that is,
the expectation with respect to the conditional distribution of
$\bm{v}$ and $\bm{\tau}$ given $\bm{\eta}^{\prime}$ and $\bm{y}$.

When calculating the upper bound of $E(\sum_{i=1}^{n}v_i^2|\bm{\eta}^{\prime}, \bm{y})$, we need to take into account the fact that $\bm{W}$ is not a full rank  matrix in the current setting. But, $E(\sum_{i=1}^{n}v_i^2|\bm{\eta}^{\prime}, \bm{y})$  can be written as,
\[
E\left[\sum_{i=1}^{n}v_{i}^{2}|\bm{\eta}^{\prime}, \bm{y}\right]=n+\sum_{i=1}^{n}\left(\tilde{\bm{w}}_{i}^{*T}\tilde{\bm{\eta}}^{\prime}\right)^{2}-\sum_{i=1}^{n}\frac{\left(\tilde{\bm{w}}_{i}^{*T}\tilde{\bm{\eta}}^{\prime}\right)\phi\left(\tilde{\bm{w}}_{i}^{*T}\tilde{\bm{\eta}}^{\prime}\right)}{1-\Phi\left(\tilde{\bm{w}}_{i}^{*T}\tilde{\bm{\eta}}^{\prime}\right)}.
\]
where $\tilde{\bm{w}}_i^*$'s are defined in section \ref{sec_propriety}. 

Since the condition B3 is in force, we know that $\tilde{\bm{W}}$ is a full rank matrix. Then the same techniques (see \eqref{eq:improper_thm2_ev}) as in the proof of Theorem \ref{them2}  can be used to show that there exists $\lambda_0 \in [0,1)$ such that 
\begin{align}
E\left[\sum_{i=1}^{n}v_{i}^{2}|\bm{\eta}^{\prime}, \bm{y}\right] & \leq\lambda_{0}\sum_{i=1}^{n}\left(\tilde{\bm{w}}_{i}^{*T}\tilde{\bm{\eta}}^{\prime}\right)^{2}+n\left(1+\Xi\right)\nonumber \\
 & =\lambda_{0}\sum_{i=1}^{n}\left(\bm{w}_{i}^{T}\bm{\eta}^{\prime}\right)^{2}+n\left(1+\Xi\right) \nonumber \\
 & =\lambda_{0}\sum_{i=1}^{n}\left(\bm{x}_{i}^{T}\bm{\beta}^{\prime}+\bm{z}_{i}^{T}\bm{u}^{\prime}\right)^{2}+n\left(1+\Xi\right). \label{eq:improper_thm4_ev}
\end{align}

For  $s\in \tilde{S}$, we have \begin{equation}
E\left[\tau_{j}^{-s}|\bm{\eta}^{\prime}, \bm{y}\right]=2^{s}G_{j}\left(s\right)\left(b_{j}+\frac{\bm{u}_{j}^{\prime T}\bm{u}_{j}^{\prime}}{2}\right)^{s}\leq G_{j}\left(s\right)\left(\bm{u}_{j}^{\prime T}\bm{u}_{j}^{\prime}\right)^{s}+2^{s}G_{j}\left(s\right)b_{j}^{s}. \label{eq:improper_thm4_taus}
\end{equation}
Also for $c\in C_1$, as in \eqref{eq:improper_thm2_tau}, we have
\begin{align*}
E\left[\tau_{j}^{c}|\bm{\eta}^{\prime}, \bm{y}\right] &=2^{-c}G_{j}\left(-c\right)\left[b_{j}+\frac{\bm{u}_{j}^{\prime T}\bm{u}_{j}^{\prime}}{2}\right]^{-c}\\
&\leq G_{j}\left(-c\right)\left[\left(2b_{j}\right)^{-c}I_{\left(0,\infty\right)}\left(b_{j}\right)+\left(\bm{u}_{j}^{\prime T}\bm{u}_{j}^{\prime}\right)^{-c}I_{\left\{ 0\right\} }\left(b_{j}\right)\right].
\end{align*}

As in the proof of Theorem \ref{them2}, we consider two cases, namely $A$ is empty and $A$ is not empty.

\textbf{Case 1}: $A$ is not empty.

Using \eqref{eq:improper_thm2_tau}, \eqref{eq:improper_thm4_ev} and
\eqref{eq:improper_thm4_taus} from \eqref{eq:improper_thm4_step1}, we
have
\begin{align}
E\left[V\left(\bm{\eta}\right)|\bm{\eta}^{\prime}\right] & =\alpha\lambda_{0}\left(1+\frac{\delta_{2}\left(s\right)}{\alpha}\right)\sum_{i=1}^{n}\left(\bm{x}_{i}^{T}\bm{\beta}^{\prime}+\bm{z}_{i}^{T}\bm{u}^{\prime}\right)^{2} \nonumber\\
 & +\delta_{3}\left(s\right)\sum_{j=1}^{r}G_{j}\left(s\right)\left(\bm{u}_{j}^{\prime T}\bm{u}_{j}^{\prime}\right)^{s}+\delta_{1}\left(c\right)\sum_{j\in A}\left(\bm{u}_{j}^{\prime T}\bm{u}_{j}^{\prime}\right)^{-c}+L_{2}\left(\alpha,s,c\right),  \label{eq:improper_thm4_step2}
\end{align}
where 
\begin{align*}
L_{2}\left(\alpha,s,c\right) & =\kappa_{1}\left(\alpha,s,c\right)+n\left(1+\Xi\right)\left(\alpha+\delta_{2}\left(s\right)\right)+\delta_{3}\left(s\right)2^{s}\sum_{j=1}^{r}G_{j}\left(s\right)b_{j}^{s}\\
+ & 2^{-c}\sum_{j\notin A}\frac{\Gamma\left(q_{j}/2-c\right)}{\Gamma\left(q_{j}/2\right)}G_{j}\left(-c\right)\left(2b_{j}\right)^{-c},
\end{align*}
and $\delta_1(c)$ is defined as in \eqref{eq:delta1c}.

We know that for $c\in C_1$, $\delta_1(c) <1$ as in Theorem \ref{them2}. Since condition 2 of Theorem \ref{them4} holds, we have $\delta_3(s) < 1$. For a fixed $s$, $\lambda_0\left( 1 + \delta_2(s)/\alpha\right) <1$ iff $\alpha > \lambda_0\delta_2(s)/(1-\lambda_0)$. So there exists a $\rho_2$ such that 
\[\rho_2  \equiv\rho_2(\alpha,s,c) =  \max\left\lbrace \lambda_0\left( 1 + \delta_2(s)/\alpha\right), \delta_3(s), \delta_1(c) \right\rbrace < 1\]
 and $L_2 \equiv L_2(\alpha,s,c) >0$ such that \eqref{drift4} holds.

\textbf{Case 2}: $A$ is empty.

In this case, the conditional expectation of $\tau_j^c$ can be bounded by a constant. Thus we have
\begin{align*}
E\left[V\left(\bm{\eta}\right)|\bm{\eta}^{\prime}\right] & =\alpha\lambda_{0}\left(1+\frac{\delta_{2}\left(s\right)}{\alpha}\right)\sum_{i=1}^{n}\left(\bm{x}_{i}^{T}\bm{\beta}^{\prime}+\bm{z}_{i}^{T}\bm{u}^{\prime}\right)^{2}\\
 & +\delta_{3}\left(s\right)\sum_{j=1}^{r}G_{j}\left(s\right)\left(\bm{u}_{j}^{\prime T}\bm{u}_{j}^{\prime}\right)^{s}+L_{2}\left(\alpha,s,c\right).
\end{align*}

As in case 1, it follows that \eqref{drift4} holds.

Since $\bm{\eta}$-chain is a Feller chain on
$\mathbb{R}^{p+q} \backslash \mathcal{N}$, and $V(\bm{\eta})$ is
unbounded off compact sets on
$\mathbb{R}^{p+q} \backslash \mathcal{N}$, the $\bm{\eta}$-chain is
geometrically ergodic on $\mathbb{R}^{p+q} \backslash
\mathcal{N}$. Using the same techniques as in Appendix
\ref{proof_thm2}, it can be shown that the original
$\{\bm{\eta}^{(m)}\}_{m=0}^{\infty}$ Markov chain defined on
$\mathbb{R}^{p+q}$ is also geometrically ergodic.
\end{proof}

\bibliographystyle{apalike} 
\bibliography{ref_probit_mixed}
\end{document}